\documentclass[11pt,reqno,a4paper]{amsart}
\usepackage[T1]{fontenc}
\usepackage{lmodern}
\usepackage{amsmath,amssymb,amsthm,mathrsfs}
\usepackage[margin=28mm]{geometry}
\usepackage{enumitem}
\usepackage[colorlinks=true,linkcolor=blue,citecolor=blue,urlcolor=blue]{hyperref}
\hypersetup{bookmarksnumbered=true,bookmarksdepth=2,pdfauthor={Yingjin Bi},
  pdftitle={Monoidal categorification of open Richardson varieties}}
\newtheorem{theorem}{Theorem}[section]
\newtheorem{proposition}[theorem]{Proposition}
\newtheorem{lemma}[theorem]{Lemma}

\theoremstyle{definition}
\newtheorem{definition}[theorem]{Definition}
\newtheorem{example}[theorem]{Example}
\newtheorem{remark}[theorem]{Remark}

\numberwithin{equation}{section}
\newcommand{\CC}{\mathbb C}

\newcommand{\ZZ}{\mathbb Z}
\newcommand{\NN}{\mathbb Z_{\geq0}}

\newcommand{\cA}{\mathcal A}

\newcommand{\cC}{\mathcal C}
\newcommand{\cS}{\mathcal S}
\newcommand{\bfa}{\mathbf a}

\newcommand{\soc}{\operatorname{soc}}

\newcommand{\Fac}{\operatorname{Fac}}
\newcommand{\Sub}{\operatorname{Sub}}

\newcommand{\Iwv}{\mathcal I_{w,v}}
\newcommand{\Kcl}[1]{K_{\mathrm{cl}}(#1)}
\newcommand{\Rwv}{\mathring{\mathcal B}_{v,w}}
\newcommand{\fr}{\mathrm{fr}}
\newcommand{\ex}{\mathrm{ex}}
\newcommand{\cl}{\mathrm{cl}}
\newcommand{\bfj}{\mathbf j}
\newcommand{\bfi}{\mathbf i}

\setlist[enumerate]{label=\textup{(\roman*)},leftmargin=*}
\setlist[itemize]{leftmargin=*}
\begin{document}
\title{Monoidal categorification of open Richardson varieties}
\author{Yingjin Bi}
\address{Department of Mathematics, Harbin Engineering University}
\email{yingjinbi@mail.bnu.edu.cn}
\date{}
\begin{abstract}
Let $G$ be a simply connected complex semisimple group of finite
simply-laced type, and let $v\leq w$ be Weyl group elements.
We construct a monoidal categorification of open Richardson
varieties in this setting. The construction
realizes the cluster variables of M\'enard's seed by real simple
modules in the quiver Hecke category $\mathscr C_{w,v}$, by
comparing their truncated Lusztig parameters with the
$\Delta$-vectors at every step of M\'enard's mutation algorithm.
Every mutation in a retained direction is realized inside
$\mathscr C_{w,v}$. After localizing at the determinantial modules
$M(w\varpi_i,v\varpi_i)$, the resulting category is a monoidal
categorification of the coordinate ring of the open Richardson
variety $\mathring{\mathcal B}_{v,w}$.
In particular, cluster monomials without negative frozen exponents
are represented by real simple modules in $\mathscr C_{w,v}$;
arbitrary frozen Laurent factors are realized in its localization.
\end{abstract}
\maketitle

\section{Introduction}

Cluster algebras connect the geometry of flag varieties with
representation theory \cite{fomin2002cluster}. A monoidal
categorification realizes a cluster algebra as the Grothendieck
ring of an abelian monoidal category, with cluster monomials
represented by real simple objects \cite{hernandez2010cluster}.
Here a simple object is called \emph{real} if its monoidal square
is simple. For quantum unipotent coordinate rings, such
categorifications are provided by categories of finite-dimensional
graded modules over symmetric quiver Hecke algebras
\cite{geiss2013cluster,kang2018monoidal}.

Let $G$ be of finite type $ADE$, and let $W$ be its Weyl group,
with longest element $w_0$ and length function $\ell$.
For $v\leq w$ in the Bruhat order, denote by $\Rwv$ the open
Richardson variety obtained by intersecting the Schubert cell
indexed by $w$ with the opposite Schubert cell indexed by $v$.
Kashiwara--Kim--Oh--Park introduced a monoidal subcategory
$\mathscr C_{w,v}\subset\mathscr C_w$, where $\mathscr C_w$ is
the category of quiver Hecke algebra modules categorifying the
quantum unipotent coordinate ring associated with $w$.
They identified the graded Grothendieck ring of
$\mathscr C_{w,v}$ with a quantum doubly invariant algebra
\cite{kashiwara2018monoidal}. After specialization at $q=1$
and localization at suitable minors, this algebra becomes
$\CC[\Rwv]$.

On the additive side, Leclerc constructed a Frobenius category
of modules over a preprojective algebra and an associated
cluster structure \cite{leclerc2016cluster}. M\'enard subsequently
gave an explicit mutation algorithm for an initial seed of this
structure \cite{menard2022cluster}. The comparison with
braid-variety seeds identifies the resulting cluster algebra,
after inversion of its frozen variables, with the full coordinate
ring $\CC[\Rwv]$ \cite[Theorem~10.1]{casals2024cluster}.
The purpose of this paper is to construct a monoidal
categorification of this cluster structure by realizing
M\'enard's mutation sequence inside $\mathscr C_w$.

The main step is a comparison between two families of integer
vectors: Lusztig parameters on the quantum unipotent side and
$\Delta$-vectors on the additive side. The former describe
dual canonical basis elements through their
Poincar\'e--Birkhoff--Witt (PBW) parametrization, whereas the
latter record filtration multiplicities of preprojective-algebra
modules. We prove that the coordinates relevant to $v$ agree
throughout M\'enard's mutation sequence. Their vanishing at
the final step then places the corresponding simple modules
in $\mathscr C_{w,v}$.

To state this comparison precisely, fix a reduced expression
$\bfi$ of $w$ and the reduced subexpression for $v$ selected
by M\'enard's algorithm, using the conventions of
Section~\ref{sec:preliminaries}. Put
$N=\ell(w_0)$ and $d=\ell(v)$, and choose a completion
$\bfj=(j_1,\ldots,j_N)$ to a reduced expression of $w_0$,
with the chosen expression of $v$ occupying the first $d$
positions. Write $\dot v=(j_N,\ldots,j_1)$ for the display
of this completed word in the additive conventions.
Thus $\dot v$ denotes the full word of length $N$.

For a quantum cluster variable $X$, let
\[
  \mathbf a^{\bfj}(X)
  =\bigl(a^{\bfj}_1(X),\ldots,a^{\bfj}_N(X)\bigr)
  \in\ZZ_{\geq0}^{N}
\]
be its \emph{Lusztig parameter} with respect to $\bfj$.
This is the exponent vector of the distinguished leading
dual PBW term in the triangular expansion of the corresponding
normalized dual canonical basis element, with the PBW
conventions fixed in Section~\ref{sec:preliminaries}.
On the additive side, for a preprojective-algebra module $R$,
write
\[
  \Delta_{\dot v}(R)
  =\bigl(\Delta_{\dot v,1}(R),\ldots,
         \Delta_{\dot v,N}(R)\bigr)
  \in\ZZ_{\geq0}^{N}
\]
for its \emph{$\Delta$-vector}. Its coordinates record the
multiplicities of the ordered layer modules in the filtration
associated with $\dot v$; these layer modules and the filtration
are recalled in Section~\ref{sec:preliminaries}.
Both vectors are indexed by the same positions
$1,\ldots,N$, with the additive positions read from right
to left in the displayed word $\dot v$.
We compare their first $d$ coordinates using the projection
\[
  \pi_v:\ZZ^N\longrightarrow\ZZ^d,
  \qquad
  (c_1,\ldots,c_N)\longmapsto(c_1,\ldots,c_d).
\]

\begin{theorem}[Comparison of truncated parameters]
\label{thm:intro-parameter-comparison}
Apply M\'enard's mutation sequence to the corresponding
additive and quantum unipotent seeds. At a step $t$, let
$R_k^{(t)}$ be the indecomposable summand of the additive
seed and $X_k^{(t)}$ the quantum cluster variable carrying
the same position label $k$ in the full seeds. Then
\[
  \pi_v\!\left(\mathbf a^{\bfj}(X_k^{(t)})\right)
  =
  \pi_v\!\left(\Delta_{\dot v}(R_k^{(t)})\right)
\]
for every step $t$ and every label $k$, including labels
omitted from the auxiliary cut seeds used in the proof.
In particular, every variable retained in the final
Richardson seed has its first $d$ Lusztig coordinates
equal to zero.
\end{theorem}

Theorem~\ref{thm:parameter-vanishing} establishes this
comparison by combining PBW leading parameters with
M\'enard's additive mutation formulas. The resulting
vanishing condition shows that the real simple modules
representing the retained variables belong to
$\mathscr C_{w,v}$. The support formula for convolution,
which is the monoidal product of quiver Hecke algebra
modules, then ensures that subsequent mutations of the
Richardson seed are also realized within this category.
Localization at determinantial modules completes the
categorification.

More precisely, let $I$ index the simple roots of $G$, and
let $\varpi_i$ be the fundamental weight corresponding to
$i\in I$. Denote by $M(w\varpi_i,v\varpi_i)$ the
determinantial simple module corresponding to the quantum
minor with extremal weights $w\varpi_i$ and $v\varpi_i$.
Write $\widetilde{\mathscr C}_{w,v}$ for the monoidal
localization of $\mathscr C_{w,v}$ in which these modules
are made invertible. Let $\Sigma_{w,v}$ be M\'enard's
initial Richardson seed, let $J_{\fr}$ be its set of frozen
indices, and write $x_a$ for the frozen variable indexed
by $a\in J_{\fr}$. We use $\cA(\Sigma_{w,v})$ for the
associated cluster algebra over $\CC$ before inversion
of the frozen variables.

For a graded monoidal category $\mathscr C$, define its
classical Grothendieck ring by
\[
  \Kcl{\mathscr C}
  :=\CC\otimes_{\ZZ[q,q^{-1}]}K_0(\mathscr C),
  \qquad q\longmapsto1,
\]
where $q$ acts on $K_0(\mathscr C)$ by grading shift.
We obtain the following theorem.

\begin{theorem}\label{thm:introduction}
For $G$ of finite type $ADE$ and $v\leq w$, the category
$\widetilde{\mathscr C}_{w,v}$ is a monoidal
categorification of $\CC[\Rwv]$ equipped with M\'enard's
cluster structure. More precisely, there are algebra
isomorphisms
\[
  \Kcl{\widetilde{\mathscr C}_{w,v}}
  \simeq
  \cA(\Sigma_{w,v})[x_a^{-1}\mid a\in J_{\fr}]
  \simeq
  \CC[\Rwv].
\]
Every cluster monomial with nonnegative frozen exponents
is the class of a real simple object of
$\mathscr C_{w,v}$. Cluster monomials with arbitrary
integer frozen exponents are represented by real simple
objects of $\widetilde{\mathscr C}_{w,v}$.
\end{theorem}

The ring isomorphisms in Theorem~\ref{thm:introduction}
are established after specialization at $q=1$ and
determinantial localization. At the quantum level, the
construction gives an inclusion of the quantum cluster
algebra with nonnegative frozen exponents into
$K_0(\mathscr C_{w,v})$, after the standard extension of
scalars allowing half-integral grading normalizations.
Equality with the unlocalized quantum Grothendieck ring
is not asserted.

Section~\ref{sec:preliminaries} fixes the conventions and
recalls the geometric, additive, and PBW constructions.
Section~\ref{sec:parameters} computes the initial
truncated parameters, and Section~\ref{sec:comparison}
compares them along M\'enard's mutation sequence.
Finally, Section~\ref{sec:categorification} constructs
the monoidal realization and proves the localization
theorem.
\subsection*{Acknowledgements}
The author thanks Fan Qin for helpful discussions.

\section{Preliminaries}
\label{sec:preliminaries}

We fix the root and word conventions, relate the quantum and
classical invariant algebras, and recall the geometric
localization and PBW parametrizations used below.

\subsection{Root and word conventions}

Throughout, $G$ is a simply connected complex semisimple group of
finite simply-laced type. Fix opposite Borel subgroups $B^+,B^-$,
their common maximal torus $H$, and their unipotent radicals $N,N^-$.
Fix a pinning, with simple root subgroups $x_i(t)\in N$ and
$y_i(t)\in N^-$. Let $I$ index the simple roots $\alpha_i$ and
coroots $\alpha_i^\vee$, and let
$C=(c_{ij})_{i,j\in I}$ be the Cartan matrix, where
$c_{ij}=\langle\alpha_i^\vee,\alpha_j\rangle$.
Write $\Delta_+$ and $\Delta_-=-\Delta_+$ for the positive and
negative roots, respectively. Set
\[
  Q_+=\sum_{i\in I}\NN\alpha_i,\qquad Q_-=-Q_+,
  \qquad P_+=\sum_{i\in I}\NN\varpi_i,
\]
where $\NN=\mathbb Z_{\geq0}$ and $\varpi_i$ is the fundamental
weight corresponding to $i$. Let $W$ be the Weyl group, with
simple reflections $s_i$, length function $\ell$, and longest
element $w_0$.

For a reduced expression $w=s_{i_1}\cdots s_{i_r}$, where
$r=\ell(w)$, write $\bfi=(i_1,\ldots,i_r)$ and use the
quiver Hecke conventions
\begin{equation}\label{eq:word-convention}
  w_{\leq k}=s_{i_1}\cdots s_{i_k},\qquad
  \beta_k=w_{\leq k-1}(\alpha_{i_k}),\qquad
  w_{\leq0}=e.
\end{equation}
Thus $\{\beta_1,\ldots,\beta_r\}=\Delta_+\cap w\Delta_-$.
The reversed tuple $\overline w=(i_r,\ldots,i_1)$ is a reduced
expression of $w^{-1}$ and is the displayed word used in the
conventions of Leclerc and M\'enard. Accordingly, our pair
$(v,w)$ corresponds to their pair $(v^{-1},w^{-1})$.
We retain the position labels $1,\ldots,r$ from the original
word, so that these labels increase from right to left in
$\overline w$.

For $k\in[1,r]$, let $k^+$ and $k^-$ be the next and previous
positions of the same color $i_k$, with boundary values
$r+1$ and $0$, respectively. Write $k^{t+}$ and $k^{t-}$
for their $t$-fold iterates, with $k^{0+}=k^{0-}=k$;
iterations remain at the boundary once it is reached.
Also put
\[
  k_{\min}=\min\{a\mid i_a=i_k\},\qquad
  k_{\max}=\max\{a\mid i_a=i_k\}.
\]
When several words occur, successors and predecessors are
computed in the word specified by the superscript or the
surrounding notation.

\subsection{Quantum coordinate rings and unipotent minors}
\label{subsec:quantum-coordinate-minors}

Put $\mathbb A=\mathbb Z[q,q^{-1}]$ and
$\mathfrak g=\operatorname{Lie}(G)$,
$\mathfrak n=\operatorname{Lie}(N)$.
Let $U_q^+(\mathfrak g)$ be the positive part of the quantum
group $U_q(\mathfrak g)$. Its restricted graded dual
\[
  A_q(\mathfrak n)
  =
  \bigoplus_{\beta\in Q_+}
  \operatorname{Hom}_{\mathbb Q(q)}
  \bigl(U_q^+(\mathfrak g)_\beta,\mathbb Q(q)\bigr)
\]
carries the standard quantum coordinate-ring multiplication.
Denote by $A_q(\mathfrak n)_{\mathbb A}$ the integral form
dual to the divided-power integral form of $U_q^+(\mathfrak g)$.
The left and right actions are specified by
\[
  (x\cdot f\cdot y)(z)=f(yzx),
  \qquad x,y,z\in U_q^+(\mathfrak g).
\]
Define
\begin{equation}\label{eq:quantum-invariant-algebra}
  A_{w,v}
  :=
  \left\{
  f\in A_q(\mathfrak n)_{\mathbb A}
  \ \middle|\
  \begin{aligned}
    U_q^+(\mathfrak g)_\beta\cdot f&=0
    &&\text{for all }\beta\in(Q_+\cap wQ_+)\setminus\{0\},\\
    f\cdot U_q^+(\mathfrak g)_\gamma&=0
    &&\text{for all }\gamma\in(Q_+\cap vQ_-)\setminus\{0\}
  \end{aligned}
  \right\}.
\end{equation}
This is the quantum doubly invariant subalgebra of
Kashiwara--Kim--Oh--Park
\cite[Section~2.2]{kashiwara2018monoidal}.

For $\lambda\in P_+$, let $V_q(\lambda)$ be the simple
highest-weight module, with normalized extremal vectors
$u_\mu$ for $\mu\in W\lambda$ and the standard contravariant
form $(\,\cdot\,,\,\cdot\,)_\lambda$.
The \emph{unipotent quantum minor} $D(\mu,\nu)$ is the
matrix coefficient
\begin{equation}\label{eq:unipotent-quantum-minor}
  D(\mu,\nu)(x)
  =(u_\nu,xu_\mu)_\lambda,
  \qquad x\in U_q^+(\mathfrak g),
  \qquad \mu,\nu\in W\lambda.
\end{equation}
The form is normalized by $(u_\lambda,u_\lambda)_\lambda=1$.
These minors belong to $A_q(\mathfrak n)_{\mathbb A}$.
Writing $\Delta_{\nu,\mu}$ for the corresponding
classical generalized minor restricted to $N$, we have
\begin{equation}\label{eq:minor-specialization}
  \left.D(\mu,\nu)\right|_{q=1}=\Delta_{\nu,\mu}.
\end{equation}
In particular, the two weight labels occur in reversed order
in the quantum and classical notation; see
\cite[Definition~1.5]{kashiwara2018monoidal}.

The upper global basis of $A_q(\mathfrak n)_{\mathbb A}$
is denoted by $\mathbf B^{\mathrm{up}}$.
A minor $D(\mu,\nu)$ is nonzero precisely when there exist
$u',v'\in W$ such that
$\mu=u'\lambda$, $\nu=v'\lambda$, and $v'\leq u'$;
in that case it belongs to $\mathbf B^{\mathrm{up}}$
\cite[Lemma~1.6]{kashiwara2018monoidal}.
Thus nonvanishing depends on the orbit weights, and
$v\leq u$ is sufficient for $D(u\lambda,v\lambda)\neq0$.
We also have $D(\mu,\mu)=1$.

\subsection{Conventions and the geometric localization}
\label{subsec:conventions-localization}

Choose representatives $\dot u\in N_G(H)$ of Weyl group
elements $u\in W$. For $v\leq w$, set
\[
  \Rwv
  :=(B^-\dot vB^+/B^+)\cap(B^+\dot wB^+/B^+).
\]
Define
\[
  N(w):=N\cap\dot wN^-\dot w^{-1},\qquad
  N'(w):=N\cap\dot wN\dot w^{-1},
\]
and specify the classical invariant ring by
\begin{equation}\label{eq:invariant-ring-convention}
  \mathcal I_{w,v}
  :=
  \{f\in\mathbb C[N]\mid
    f(anb)=f(n)
    \text{ for all }a\in N(v),\ b\in N'(w),\ n\in N\}.
\end{equation}
With the actions used in
\eqref{eq:quantum-invariant-algebra}, specialization gives
\[
  \mathbb C\otimes_{\mathbb A}A_{w,v}
  \simeq\mathcal I_{w,v},
  \qquad q\longmapsto1;
\]
see \cite[Remark~2.19 and Theorem~2.20]{kashiwara2018monoidal}.

Let $\Lambda$ be the preprojective algebra of the Dynkin diagram,
$S_i$ its simple module at $i$, and $Q_i$ the injective envelope
of $S_i$. For a sequence $(j_1,\ldots,j_t)$, the iterated socle
$\soc_{(j_1,\ldots,j_t)}(M)$ is constructed by successively
taking the $S_{j_a}$-isotypical socle in the quotient by the
submodule already obtained. The functor $\mathcal E_i$
removes the $S_i$-isotypical head, while
$\mathcal E_i^\dagger$ takes the quotient by the
$S_i$-isotypical socle. Their braid relations give well-defined
functors
\[
  \mathcal E_u=\mathcal E_{a_1}\cdots\mathcal E_{a_t},
  \qquad
  \mathcal E_u^\dagger
  =\mathcal E_{a_1}^\dagger\cdots\mathcal E_{a_t}^\dagger
\]
for a reduced expression $u=s_{a_1}\cdots s_{a_t}$;
see \cite[Section~3.2]{leclerc2016cluster}.

Leclerc's categories are
\[
  \mathcal C_b
  =\Fac\!\left(\mathcal E_{b^{-1}w_0}
      \left(\bigoplus_{i\in I}Q_i\right)\right),
  \qquad
  \mathcal C^a
  =\Sub\!\left(\mathcal E_{a^{-1}}^\dagger
      \left(\bigoplus_{i\in I}Q_i\right)\right),
\]
where $\Fac$ and $\Sub$ denote factor modules and submodules
of finite direct sums, respectively. Put
$\mathcal C^{\mathrm{Lec}}_{a,b}=\mathcal C_b\cap\mathcal C^a$.
In our conventions, the relevant additive category is
\[
  \mathcal C^{\mathrm{add}}_{w,v}
  :=\mathcal C^{\mathrm{Lec}}_{v^{-1},w^{-1}}.
\]
This is a Frobenius category with a cluster structure
\cite[Corollary~3.11]{leclerc2016cluster}.
The rightmost reduced subexpression for $v^{-1}$ in
$(i_r,\ldots,i_1)$ corresponds to the leftmost reduced
subexpression for $v$ in $(i_1,\ldots,i_r)$, with the same
position labels $p_1<\cdots<p_{\ell(v)}$.

Let $g\mapsto g^{\mathsf T}$ be the transpose
anti-automorphism determined by
\[
  x_i(t)^{\mathsf T}=y_i(t),\qquad
  y_i(t)^{\mathsf T}=x_i(t),\qquad
  h^{\mathsf T}=h\quad(h\in H).
\]
It induces an isomorphism
\begin{equation}\label{eq:left-right-richardson}
  \mathcal R^{\mathrm{Lec}}_{v^{-1},w^{-1}}
  \xrightarrow{\ \sim\ }\Rwv,
  \qquad B^-g\longmapsto g^{\mathsf T}B^+.
\end{equation}
Indeed, transpose exchanges $B^-$ and $B^+$, reverses
products, and sends a representative of $u$ to a
representative of $u^{-1}$. Moreover,
$N_{\mathrm{Lec}}(v^{-1})=N(v)$ and
$N'_{\mathrm{Lec}}(w^{-1})=N'(w)$.
Hence Leclerc's invariant ring for $(v^{-1},w^{-1})$
is precisely $\mathcal I_{w,v}$, with no additional
involution on $\mathbb C[N]$.

Let $\varphi$ be the cluster character for
preprojective-algebra modules. For
\[
  V_k=\soc_{(i_k,\ldots,i_1)}(Q_{i_k})
  \in\mathcal C_{w^{-1}},
\]
its normalization is
\[
  \varphi(V_k)
  =\Delta_{\varpi_{i_k},w_{\leq k}\varpi_{i_k}}
  =\left.D(w_{\leq k}\varpi_{i_k},\varpi_{i_k})
    \right|_{q=1}.
\]
Thus the common position labels identify the additive
initial cluster variables with the classical
specializations of the quantum unipotent initial
variables. Applying the same mutation sequence preserves
this identification.

We next explain why the frozen variables are independent
of the categorical seed. Let $T$ be a basic cluster-tilting
object of $\mathcal C^{\mathrm{add}}_{w,v}$, and let $P$
be projective-injective in this category. Since
$\operatorname{Ext}^1_\Lambda(T,P)=0$, the cluster-tilting
property implies $P\in\operatorname{add}(T)$.
Here $\operatorname{add}(T)$ consists of direct summands
of finite direct sums of copies of $T$.
Consequently, every indecomposable projective-injective
occurs, up to isomorphism, exactly once in every basic
cluster-tilting object. Categorical mutation replaces only
non-projective-injective summands, so these common
summands determine the frozen variables in every
categorical seed of $\mathcal C^{\mathrm{add}}_{w,v}$.

Choose representatives $P_a$, $a\in J_{\mathrm{fr}}$, of
these indecomposable projective-injectives, and put
$x_a=\varphi(P_a)\in\mathcal I_{w,v}$.
The inclusion follows from
\cite[Proposition~4.2]{leclerc2016cluster}.
The following proposition relates this common frozen set
to the minors defining the geometric localization.

\begin{proposition}\label{prop:geometric-frozen-localization}
For $i\in I$, set
\[
  P_i^{\mathrm{det}}
  :=\mathcal E_v^\dagger\mathcal E_{ww_0}(Q_i),
  \qquad
  d_i:=\Delta_{v\varpi_i,w\varpi_i}
      =\left.D(w\varpi_i,v\varpi_i)\right|_{q=1}.
\]
Then
$\operatorname{add}(\bigoplus_{i\in I}P_i^{\mathrm{det}})$
is the subcategory of projective-injective objects of
$\mathcal C^{\mathrm{add}}_{w,v}$, and
$\varphi(P_i^{\mathrm{det}})=d_i$.
Writing
\[
  P_i^{\mathrm{det}}
  \simeq\bigoplus_{a\in J_{\mathrm{fr}}}
          P_a^{\oplus m_{ia}},
  \qquad m_{ia}\in\mathbb Z_{\geq0},
\]
we obtain
\[
  d_i=\prod_{a\in J_{\mathrm{fr}}}x_a^{m_{ia}}.
\]
Every $x_a$ occurs with positive exponent in at least one
of these products. In particular,
\[
  \mathcal I_{w,v}[d_i^{-1}\mid i\in I]
  =
  \mathcal I_{w,v}[x_a^{-1}\mid a\in J_{\mathrm{fr}}]
  \simeq\mathbb C[\Rwv].
\]
\end{proposition}

\begin{proof}
By \cite[Lemma~3.15 and Proposition~4.1]{leclerc2016cluster},
applied to $(v^{-1},w^{-1})$, the objects
$\mathcal E_v^\dagger\mathcal E_{ww_0}(Q_i)$ generate
the projective-injective subcategory, and their cluster
characters are $\Delta_{v\varpi_i,w\varpi_i}$.
Their indecomposable summands are therefore precisely
the objects $P_a$ that occur as frozen summands in every
categorical seed.

The stated factorizations now follow from
$\varphi(M\oplus M')=\varphi(M)\varphi(M')$.
Since the $P_i^{\mathrm{det}}$ generate the
projective-injective subcategory, for every $a$ there is
some $i$ with $m_{ia}>0$. For such an $i$,
\[
  x_a^{-1}
  =
  d_i^{-1}x_a^{m_{ia}-1}
  \prod_{\substack{b\in J_{\mathrm{fr}}\\b\neq a}}
       x_b^{m_{ib}}.
\]
Hence inverting the $d_i$ makes every $x_a$ invertible.
Conversely, inverting all $x_a$ makes every $d_i$
invertible, proving equality of the two localizations.

Finally, the geometric localization theorem
\cite[Theorem~2.12]{leclerc2016cluster}, together with
\eqref{eq:left-right-richardson}, identifies
$\mathcal I_{w,v}[d_i^{-1}\mid i\in I]$
with $\mathbb C[\Rwv]$.
\end{proof}

Thus localization at the determinantial minors agrees
with localization at the common frozen variables of the
categorical seeds. This statement concerns the invariant
ring and its localization; equality between
$\mathcal I_{w,v}$ and the cluster algebra before
inversion of frozen variables is not asserted.

\subsection{Interval minors and PBW parameters}
\label{subsec:interval-minors-pbw}

We now describe the minors attached to subintervals of
$\bfi$ and the PBW parameters used to compare additive
and monoidal mutations. All quantum minors use the
normalization of \eqref{eq:unipotent-quantum-minor}.

An \emph{$i$-box} in $\bfi$ is an interval $[a,b]$ with
$1\leq a\leq b\leq r$ and $i_a=i_b=i$. Set
\begin{equation}\label{eq:interval-minor}
  D^{\bfi}[a,b]
  :=D(w_{\leq b}\varpi_i,w_{\leq a-1}\varpi_i).
\end{equation}
Since all letters strictly between $a^-$ and $a$ have
color different from $i$, their simple reflections fix
$\varpi_i$. Hence
$w_{\leq a-1}\varpi_i=w_{\leq a^-}\varpi_i$,
where $w_{\leq0}=e$ if $a^-=0$.
For an empty interval, we set $D^{\bfi}[a,b]=1$.
In particular, the initial minors are
\[
  D_k:=D(w_{\leq k}\varpi_{i_k},\varpi_{i_k})
      =D^{\bfi}[k_{\min},k],
  \qquad 1\leq k\leq r.
\]
A singleton interval $D^{\bfi}[a,a]$ is the dual root
vector associated with $\beta_a$ in the normalization
used below; it is not the empty-interval minor.

Put $L=\ell(w_0)$ and choose a completion
$\widehat{\bfi}=(i_1,\ldots,i_L)$ of $\bfi$ to a reduced
expression of $w_0$. Extend the definition of $\beta_t$
in \eqref{eq:word-convention} to $1\leq t\leq L$.
We suppress the hat in the PBW notation associated with
this fixed completion. Let $E_{\bfi}^*(m\beta_t)$ denote
the normalized dual divided-power root vector, with
$E_{\bfi}^*(0)=1$ and
$E_{\bfi}^*(\beta_t)=D^{\widehat{\bfi}}[t,t]$.

We use the decreasing-root PBW convention, corresponding
to $e=-1$ in \cite[Section~4.3]{kimura2012quantum}.
In particular, both the root vectors and the product
order are taken in this convention.
For $\mathbf a=(a_1,\ldots,a_L)\in\NN^L$, write
\[
  P_{\bfi}(\mathbf a)
  =
  q^{\kappa_{\bfi}(\mathbf a)}
  E_{\bfi}^*(a_L\beta_L)\cdots
  E_{\bfi}^*(a_1\beta_1).
\]
Here $\kappa_{\bfi}(\mathbf a)$ is the scalar
normalization exponent, chosen so that the transition
from the upper global basis to this dual PBW basis has
diagonal coefficients equal to $1$.
Let $G_{\bfi}^{\mathrm{up}}(\mathbf a)$ be the upper
global basis element with Lusztig parameter $\mathbf a$.

We order $\NN^L$ lexicographically in increasing
coordinate order: $\mathbf b<_{\mathrm{lex}}\mathbf a$
if, at the first coordinate where they differ, the
entry of $\mathbf b$ is smaller.
With the preceding conventions, PBW triangularity gives
\begin{equation}\label{eq:pbw-global-triangularity}
  G_{\bfi}^{\mathrm{up}}(\mathbf a)
  =
  P_{\bfi}(\mathbf a)
  +
  \sum_{\mathbf b<_{\mathrm{lex}}\mathbf a}
    c_{\mathbf b,\mathbf a}(q)P_{\bfi}(\mathbf b).
\end{equation}
The sum is finite, and the coefficients belong to
$\mathbb A$. The triangularity and multiplication
statements of
\cite[Theorem~4.26 and Proposition~4.28]{kimura2012quantum}
apply in this decreasing-root convention by transport
through the standard anti-involution.

For a nonzero element with PBW expansion
$f=\sum_{\mathbf a}f_{\mathbf a}(q)P_{\bfi}(\mathbf a)$,
define its \emph{leading PBW parameter} by
\[
  \mathbf a^{\bfi}(f)
  :=
  \max_{\mathrm{lex}}
  \{\mathbf a\mid f_{\mathbf a}(q)\neq0\}.
\]
For an upper global basis element, this is its Lusztig
parameter by \eqref{eq:pbw-global-triangularity}.
The PBW multiplication formula yields
\begin{equation}\label{eq:pbw-leading-product}
  \mathbf a^{\bfi}(fg)
  =
  \mathbf a^{\bfi}(f)+\mathbf a^{\bfi}(g)
  \qquad(f,g\neq0).
\end{equation}
Indeed, PBW straightening is triangular for an order
compatible with addition, and the product of the leading
terms has a nonzero leading coefficient.
When quantum seed normalizations require half-integral
powers of $q$, we extend scalars to $\mathbb Q(q^{1/2})$.
The same definitions and formulas apply, and multiplying
an element by a power of $q^{1/2}$ does not change its
leading parameter. We use the same notation
$\mathbf a^{\bfj}$ for any other reduced expression
$\bfj$ of $w_0$.

For the interval minor in \eqref{eq:interval-minor},
\cite[Proposition~4.6 and Corollary~4.7]{kashiwara2018monoidal}
give
\begin{equation}\label{eq:interval-parameter}
  a_t^{\bfi}(D^{\bfi}[a,b])
  =
  \begin{cases}
    1,&a\leq t\leq b\text{ and }i_t=i_a,\\
    0,&\text{otherwise},
  \end{cases}
  \qquad 1\leq t\leq L.
\end{equation}
In particular, all coordinates beyond $r$ vanish.
Thus the parameter may be regarded as a vector in
$\NN^r$, independent of the chosen completion, and then
extended by zero to $\NN^L$. For $a=b$, it is the standard
unit vector at position $a$; for an empty interval, it
is the zero vector.

We will use the following consequence of
\eqref{eq:pbw-leading-product} for quantum mutations.

\begin{lemma}[Leading parameters in an exchange relation]
\label{lem:leading-exchange}
Let $(X_i)_{i\in J}$ be the variables of a quantum seed
obtained by mutations from a quantum unipotent seed,
with frozen variables not inverted, and let
$\widetilde B=(b_{ik})$ be its extended exchange matrix.
Fix a reduced expression $\bfj$ of $w_0$.
If $k$ is mutable and $X_k'$ is the variable obtained by
mutation at $k$, then
\begin{equation}\label{eq:leading-exchange}
  \mathbf a^{\bfj}(X_k)+\mathbf a^{\bfj}(X_k')
  =
  \max_{\mathrm{lex}}
  \left\{
    \sum_{i\in J}[b_{ik}]_+\mathbf a^{\bfj}(X_i),
    \quad
    \sum_{i\in J}[-b_{ik}]_+\mathbf a^{\bfj}(X_i)
  \right\},
\end{equation}
where $[m]_+=\max\{m,0\}$.
The maximum is taken in the specified lexicographic
order, not coordinatewise.
\end{lemma}

\begin{proof}
By the monoidal categorification of quantum unipotent
coordinate rings \cite{kang2018monoidal}, the quantum
cluster variables and normalized cluster monomials
belong to $\mathbf B^{\mathrm{up}}$ up to the prescribed
powers of $q^{1/2}$. Their leading PBW coefficients are
therefore powers of $q^{1/2}$ with positive sign.

Write the exchange relation as
\[
  X_kX_k'=q^{s_+}M_++q^{s_-}M_-,
  \qquad s_+,s_-\in\tfrac12\mathbb Z,
\]
where $M_+$ and $M_-$ are the normalized cluster
monomials with exponent vectors
$([b_{ik}]_+)_{i\in J}$ and
$([-b_{ik}]_+)_{i\in J}$, respectively.
By \eqref{eq:pbw-leading-product}, their leading
parameters are the two sums on the right-hand side
of \eqref{eq:leading-exchange}.

If these parameters differ, the larger one is the
leading parameter of the exchange sum. If they agree,
the leading coefficients add without cancellation,
since each is a Laurent monomial in $q^{1/2}$ with
coefficient $1$.
Thus the leading parameter of the right-hand side
is their lexicographic maximum.
Applying \eqref{eq:pbw-leading-product} to $X_kX_k'$
proves the formula.
\end{proof}

\subsection{Seeds and frozen variables}
An ice seed has a vertex set $J=J_{\ex}\sqcup J_{\fr}$ and an
extended exchange matrix $B=(b_{ij})_{i\in J,j\in J_{\ex}}$.
For a drawn quiver we use
$b_{ij}=\#\{i\to j\}-\#\{j\to i\}$.
Throughout, $\cA(\Sigma)$ denotes the cluster algebra with
\emph{noninvertible} frozen variables. We distinguish it from
\[
 \cA^{\mathrm{loc}}(\Sigma)
 :=\cA(\Sigma)[x_a^{-1}\mid a\in J_{\fr}].
\]
The same distinction is made for quantum seeds. Classical cluster
algebras are taken over $\CC$; quantum ones use the standard
$\ZZ[q^{\pm1/2}]$ coefficient extension when needed.

The ambient unipotent seed $\cS$ associated with $\bfi$ has variables
$D_1,\ldots,D_r$ and frozen positions
\[
 J_{\fr}^{\mathrm{amb}}=\{k\in[1,r]\mid k^+=r+1\}.
\]
Its exchange matrix and compatible commutation matrix are those
of the unipotent seed in
\cite{geiss2013cluster,kang2018monoidal}. One may draw its quiver
with horizontal arrows $k\to k^-$ and ordinary arrows
$s\to t$ for
\[
 s<t<s^+\leq t^+,\qquad c_{i_s,i_t}=-1,
\]
using the corresponding sign convention for $B$.
This is the opposite of M\'enard's endomorphism quiver. A
simultaneous reversal of all arrows preserves the mutation labels
and the unordered pair of exchange monomials.
Arrows between two frozen vertices do not enter the extended
exchange matrix.

Under $q=1$, the initial quantum minor $D_k$ becomes
$\Delta_{\varpi_{i_k},w_{\leq k}\varpi_{i_k}}|_N$.
The cluster algebra associated with this classical seed is
$\CC[N(w)]$.
On the additive side, the same function is the GLS cluster
character of the module associated with the reversed tuple
$\overline w$; this identification will be used in
Section~\ref{sec:comparison}.

\begin{lemma}[Freezing and deletion]\label{lem:freezing-deletion}
Let $E$ be a set of vertices to be deleted from a seed. Freeze
every retained mutable vertex that is joined to a vertex of $E$,
and retain every original frozen vertex not in $E$.
For any remaining mutable vertex $k$, the entries $b_{ek}$ vanish
for all $e\in E$. This vanishing is preserved by every sequence of
mutations in the remaining mutable directions. Such mutations
therefore agree with the corresponding ambient mutations after
deletion of $E$.
\end{lemma}
\begin{proof}
The initial vanishing is exactly the freezing rule. If $k$ is
mutated and $j$ is another retained mutable vertex, the matrix
mutation formula changes $b_{ej}$ by an expression containing
$b_{ek}=0$. It also changes $b_{ek}$ to its negative, which is still
zero. Induction proves the assertion. The exchange relations then
involve only retained variables at every step.
\end{proof}

\section{Lusztig parameters and the Richardson seed}
\label{sec:parameters}

\subsection{Lusztig parameters of determinantial minors}
\label{sec:lusD0k}

Throughout this section, we work in finite simply-laced type and retain the word convention in
\eqref{eq:word-convention}. For compatibility with the reversed
displayed words, we also write
\[
 D^{\overline w}\{0,k]:=D_k,\qquad
 \mathbf a^{\overline w}:=\mathbf a^{\mathbf i}.
\]
Thus the displayed word is reversed, but the coordinates are
always read in increasing position order.

Put $n_0:=\ell(w_0)$. We use stable MV polytopes, that is,
MV polytopes modulo translation, in the root realization
compatible with our PBW convention. Each translation class is
normalized by placing its lowest vertex at $0$.
A reduced word of $w_0$ determines a path whose successive
edge directions are the corresponding ordered positive roots.
The edge lengths, with zero lengths allowed, form its Lusztig
datum.

\begin{theorem}[{\cite[Theorem~7.2]{kamnitzer2010mirkovic}}]
\label{them:MVpolytope}
The crystal labels of the upper global basis are parametrized
by stable MV polytopes. For every reduced word of $w_0$, the
Lusztig parameters of a basis element are the edge lengths
along the corresponding path in its normalized MV polytope.
\end{theorem}

Choose a completion
\[
    \dot w=(i_{n_0},\ldots,i_{r+1},i_r,\ldots,i_1),
    \qquad s_{i_1}\cdots s_{i_{n_0}}=w_0.
\]
For a fixed $1\leq k\leq r$, put
\[
    u_k:=s_{i_{n_0}}\cdots s_{i_{k+1}}=w_0w_{\leq k}.
\]
Let $\dot x=(j_{n_0},\ldots,j_1)$ be another reduced word of $w_0$,
again displayed in reverse order. We write
$\mathbf a^{\dot x}:=\mathbf a^{(j_1,\ldots,j_{n_0})}$;
thus coordinates are indexed by the forward order.
Write
\[
    1\leq q_1<\cdots<q_{{n_0}-k}\leq {n_0}
\]
for the positions of the leftmost reduced subexpression of $u_k$
in this displayed word, so that
\[
    u_k=s_{j_{q_{{n_0}-k}}}\cdots s_{j_{q_1}}.
\]
Denote its complementary positions by
\[
    \{r_1<\cdots<r_k\}
    =\{1,\ldots,{n_0}\}\setminus\{q_1,\ldots,q_{{n_0}-k}\}.
\]
If $k={n_0}$, the $q$-sequence is empty and the complementary sequence
is $(1,\ldots,{n_0})$.

\begin{definition}\label{def:extremal-weight-recursion}
For the fixed minor index $k$, define
$\xi_t^{\dot x,k}$, $0\leq t\leq {n_0}$, by
\begin{equation}\label{eq:xi}
    \xi_0^{\dot x,k}:=\varpi_{i_k},\qquad
    \xi_t^{\dot x,k}:=
    \begin{cases}
        s_{\beta_t^{\dot x}}\xi_{t-1}^{\dot x,k},
            &t\in\{r_1,\ldots,r_k\},\\
        \xi_{t-1}^{\dot x,k},&t\in\{q_1,\ldots,q_{{n_0}-k}\},
    \end{cases}
\end{equation}
where
$\beta_t^{\dot x}=s_{j_1}\cdots s_{j_{t-1}}(\alpha_{j_t})$.
\end{definition}

The following is the extremal-minor form of the MV-polytope
calculation in \cite[\S4.1]{naito2009mirkovic}; the corresponding
preprojective-module calculation and its identification with
Lusztig data are described in
\cite[Example~5.14, Proposition~5.24 and Remark~5.25(ii)]
{baumann2014affine}.  See also
\cite[Definition~5.20 and Proposition~5.21]{menard2022cluster}
for the reverse-word notation used here.

\begin{proposition}\label{prop:extremal-minor-Lusztig-data}
For $1\leq t\leq {n_0}$, define $n_t^{\dot x,k}$ by
\[
    \xi_{t-1}^{\dot x,k}-\xi_t^{\dot x,k}
    =n_t^{\dot x,k}\beta_t^{\dot x}.
\]
Then $n_t^{\dot x,k}\in\mathbb Z_{\geq0}$ and
\[
    a_t^{\dot x}(D^{\overline w}\{0,k])
    =n_t^{\dot x,k}.
\]
\end{proposition}

\begin{proof}
Put $N=n_0$ and $\lambda=\varpi_{i_k}$, and write
\[
  D_k:=D^{\overline w}\{0,k]
      =D(w_{\leq k}\lambda,\lambda).
\]
All positions below are numbered in the forward order
$j_1,\ldots,j_N$, although $\dot x$ is displayed in reverse
order.

\smallskip
\noindent\emph{Step 1: a formula for the recursively defined weights.}
Let $\mathcal Q=\{q_1,\ldots,q_{N-k}\}$ and put
\[
  h_t:=s_{j_1}\cdots s_{j_t},\qquad
  z_t:=\prod_{\substack{q\in\mathcal Q\\q\leq t}}s_{j_q},
  \qquad h_0=z_0=e,
\]
where the factors in $z_t$ are taken in increasing order of $q$.
Then $z_N=u_k^{-1}$.
The leftmost reduced subexpression in the reverse display
is the rightmost reduced subexpression in the forward word.
Consequently, its partial products are obtained by the
backward greedy rule
\begin{equation}\label{eq:extremal-backward-greedy}
  z_{t-1}
  =
  \begin{cases}
    z_ts_{j_t},&\ell(z_ts_{j_t})<\ell(z_t),\\
    z_t,&\ell(z_ts_{j_t})>\ell(z_t).
  \end{cases}
\end{equation}
The first case occurs precisely when $t\in\mathcal Q$.

We claim that
\begin{equation}\label{eq:extremal-weight-closed-form}
  \xi_t^{\dot x,k}=h_tz_t^{-1}\lambda
  \qquad(0\leq t\leq N).
\end{equation}
This is clear for $t=0$. If $t\in\mathcal Q$, then
$z_t=z_{t-1}s_{j_t}$, and hence
\[
  h_tz_t^{-1}
  =h_{t-1}s_{j_t}s_{j_t}z_{t-1}^{-1}
  =h_{t-1}z_{t-1}^{-1}.
\]
If $t\notin\mathcal Q$, then $z_t=z_{t-1}$ and
\[
  h_tz_t^{-1}
  =
  \bigl(h_{t-1}s_{j_t}h_{t-1}^{-1}\bigr)
  h_{t-1}z_{t-1}^{-1}
  =
  s_{\beta_t^{\dot x}}h_{t-1}z_{t-1}^{-1}.
\]
These are exactly the two cases of \eqref{eq:xi}, proving
\eqref{eq:extremal-weight-closed-form}.

\smallskip
\noindent\emph{Step 2: nonnegativity.}
The reflection formula and
\eqref{eq:extremal-weight-closed-form} give
\begin{equation}\label{eq:extremal-coefficient-pairing}
  n_t^{\dot x,k}
  =
  \begin{cases}
    \langle z_{t-1}\alpha_{j_t}^{\vee},\lambda\rangle,
      &t\notin\mathcal Q,\\
    0,&t\in\mathcal Q.
  \end{cases}
\end{equation}
Indeed, in the first case,
\[
  \langle(\beta_t^{\dot x})^\vee,\xi_{t-1}^{\dot x,k}\rangle
  =
  \langle h_{t-1}\alpha_{j_t}^\vee,
          h_{t-1}z_{t-1}^{-1}\lambda\rangle
  =
  \langle z_{t-1}\alpha_{j_t}^\vee,\lambda\rangle.
\]
When $t\notin\mathcal Q$, the greedy rule says that
$z_{t-1}=z_t$ and $\ell(z_ts_{j_t})>\ell(z_t)$.
Thus $z_{t-1}\alpha_{j_t}$ is a positive root.
Since $\lambda$ is dominant integral,
\eqref{eq:extremal-coefficient-pairing} proves
$n_t^{\dot x,k}\in\mathbb Z_{\geq0}$.

\smallskip
\noindent\emph{Step 3: compatibility with braid moves.}
Consider a three-term braid move in the forward word,
replacing the letters in positions $m+1,m+2,m+3$ by
\[
  (i,j,i)\longleftrightarrow(j,i,j),
  \qquad c_{ij}=c_{ji}=-1.
\]
The backward greedy algorithm processes the same suffix
in both words, so its residual $z_{m+3}$ is the same.
Take the right parabolic decomposition
\[
  z_{m+3}=zy,\qquad
  y\in W_{\{i,j\}},\qquad
  z\alpha_i,z\alpha_j\in\Delta_+.
\]
Here $z$ is the minimal representative of its right
$W_{\{i,j\}}$-coset, so
$\ell(zy')=\ell(z)+\ell(y')$ for every
$y'\in W_{\{i,j\}}$.
The greedy decisions inside the braid interval therefore
depend only on $y$. Processing either reduced expression
of the longest element of $W_{\{i,j\}}$ removes $y$ and
leaves $z$. In particular, $z_m=z$ for both words, and
the selected positions outside the interval agree.

Set
\[
  A:=\langle z\alpha_i^\vee,\lambda\rangle,\qquad
  B:=\langle z\alpha_j^\vee,\lambda\rangle.
\]
Both numbers are nonnegative. Write $1$ for a selected
position in $\mathcal Q$ and $0$ for a complementary
position. The six possible elements $y\in W_{\{i,j\}}$
give the following exhaustive table. The coefficient
triples are computed from
\eqref{eq:extremal-coefficient-pairing}, using
$s_i\alpha_j=\alpha_i+\alpha_j$ and
$s_j\alpha_i=\alpha_i+\alpha_j$:
\[
\begin{array}{c|c|c|c}
 \text{selection in }iji&
 \text{selection in }jij&
 \text{coefficients for }iji&
 \text{coefficients for }jij\\ \hline
 000&000&(A,B,A)&(B,A,B)\\
 001&010&(A,B,0)&(B,0,A+B)\\
 010&001&(A,0,A+B)&(B,A,0)\\
 011&110&(A,0,0)&(0,0,A)\\
 110&011&(0,0,B)&(B,0,0)\\
 111&111&(0,0,0)&(0,0,0)
\end{array}
\]
For example, the selection $010$ gives successive
selected partial products
$z,z,zs_j,zs_j$, and hence the coefficients
$(A,0,A+B)$.

Since $A,B\geq0$, every row satisfies
\begin{equation}\label{eq:extremal-three-move}
  (a,b,c)\longmapsto
  (b+c-p,p,a+b-p),
  \qquad p=\min\{a,c\}.
\end{equation}
For $t\leq m$ and $t\geq m+3$, both $h_t$ and $z_t$
agree for the two words. Hence
\eqref{eq:extremal-weight-closed-form} shows that the weights
at these positions agree, and all coefficients outside the
braid interval are unchanged.

For a commuting move, the same argument applies to the
rank-two parabolic subgroup of type $A_1\times A_1$.
The two roots are orthogonal, so the affected
coefficients are interchanged and all other coefficients
remain unchanged.

Thus the vectors $(n_t^{\dot x,k})_{t=1}^N$ transform
under every braid move by the usual Lusztig transition
maps; see
\cite[Propositions~5.26 and~5.27]{baumann2014affine}.

\smallskip
\noindent\emph{Step 4: identification with the minor parameters.}
For the original completion $\dot w$, the chosen leftmost
subexpression of $u_k$ occupies positions $k+1,\ldots,N$.
Consequently,
\[
  \xi_t^{\dot w,k}
  =
  \begin{cases}
    w_{\leq t}\lambda,&0\leq t\leq k,\\
    w_{\leq k}\lambda,&k\leq t\leq N.
  \end{cases}
\]
It follows that
\[
  n_t^{\dot w,k}
  =
  \begin{cases}
    \langle\alpha_{i_t}^\vee,\lambda\rangle
      =\delta_{i_t,i_k},&t\leq k,\\
    0,&t>k.
  \end{cases}
\]
By \eqref{eq:interval-parameter}, these are precisely
the Lusztig parameters of $D_k$ in the word $\dot w$.

Any two reduced expressions of $w_0$ are connected by
commuting and three-term braid moves.
The recursion coefficients and the Lusztig parameters
of the fixed minor $D_k$ agree for $\dot w$ and undergo
the same transition at each such move.
They therefore agree for every $\dot x$, proving
\[
  a_t^{\dot x}(D^{\overline w}\{0,k])
  =n_t^{\dot x,k}
  \qquad(1\leq t\leq N).
\]
\end{proof}

The corresponding normalized MV path has vertices
\[
  \mu_t=\varpi_{i_k}-\xi_t^{\dot x,k},
  \qquad
  \mu_t-\mu_{t-1}=n_t^{\dot x,k}\beta_t^{\dot x}.
\]
In particular, $\mu_0=0$ and
$\mu_{n_0}=\varpi_{i_k}-w_{\leq k}\varpi_{i_k}$.
This explains the sign in the weight recursion.

\begin{example}\label{exam:minor-Lusztig-A3}
In type $A_3$, take
\[
    \overline w_0=(2,1,2,3,2,1),\qquad
    \dot x=(1,2,3,1,2,1),\qquad k=3.
\]
Then $u_3=s_2s_1s_2=s_1s_2s_1$.  Its leftmost reduced
subexpression in $\dot x$ has positions
\[
    (q_3,q_2,q_1)=(6,5,3),\qquad
    (r_1,r_2,r_3)=(1,2,4).
\]
The root sequence associated with $\dot x$ is
\[
    (\beta_1^{\dot x},\ldots,\beta_6^{\dot x})
    =(\alpha_1,\alpha_1+\alpha_2,\alpha_2,
      \alpha_1+\alpha_2+\alpha_3,\alpha_2+\alpha_3,\alpha_3).
\]
Since $i_3=3$, the recursion gives
\[
    \xi_0^{\dot x,3}=\xi_1^{\dot x,3}
    =\xi_2^{\dot x,3}=\xi_3^{\dot x,3}=\varpi_3,
\]
and
\[
    \xi_4^{\dot x,3}=\xi_5^{\dot x,3}=\xi_6^{\dot x,3}
    =\varpi_3-\alpha_1-\alpha_2-\alpha_3.
\]
Consequently,
\[
    \mathbf a^{\dot x}(D^{\overline w_0}\{0,3])
    =(0,0,0,1,0,0).
\]
This can also be checked by the braid moves
$(1,2,3,2,1,2)\mapsto(1,2,3,1,2,1)
\mapsto(1,2,1,3,2,1)$ on the words written in increasing
position order: the original parameter $(0,0,1,0,0,0)$
is carried to $(0,0,0,1,0,0)$.
\end{example}

\subsection{The selected subexpression and its notation}

Fix $v\leq w$ and write $w=s_{i_1}\cdots s_{i_r}$.
Put $d:=\ell(v)$ and retain the reverse display
$\overline w=(i_r,\ldots,i_1)$.
Let
\[
  \overline v=(i_{p_d},\ldots,i_{p_1}),\qquad
  1\leq p_1<\cdots<p_d\leq r,
\]
be the rightmost reduced subexpression of $v^{-1}$ in
$\overline w$. Equivalently,
$(i_{p_1},\ldots,i_{p_d})$ is the leftmost reduced
subexpression of $v$ in $(i_1,\ldots,i_r)$, and
\[
  v=s_{i_{p_1}}\cdots s_{i_{p_d}}.
\]
Thus M\'enard's Weyl group elements are
$v_{\mathrm M}=v^{-1}$ and $w_{\mathrm M}=w^{-1}$.

For $1\leq k\leq r$ and $0\leq m\leq d$, define
\begin{align*}
  f_{\min}(k)
  &:=\min\bigl(\{1\leq j\leq d\mid i_{p_j}=i_k\}
                    \cup\{d+1\}\bigr),\\
  f(k)
  &:=\max\bigl(\{1\leq j\leq d\mid p_j\leq k,\
                         i_{p_j}=i_k\}\cup\{0\}\bigr),\\
  \alpha(k,m)
  &:=\#\{1\leq j\leq m\mid i_{p_j}=i_k\}.
\end{align*}
For $1\leq m\leq d$, put
\[
  \gamma_m:=\alpha(p_m,m),\qquad
  b_m:=\#\{1\leq k\leq p_m\mid
      i_k=i_{p_m},\ k\notin\{p_1,\ldots,p_m\}\}.
\]
Thus $\gamma_m$ counts the selected positions of color
$i_{p_m}$ up to $p_m$, whereas $b_m$ counts the unselected
positions of that color preceding $p_m$.
We write $b_m$ for M\'enard's integer $\beta_m$ to
distinguish it from the positive roots.

Successors, predecessors, and their iterates in the original
word retain the notation of \eqref{eq:word-convention}.
The boundary values are fixed under further iteration.
For the selected word, write $j^\oplus$ for the next index
of the same color in the increasing order $1,\ldots,d$,
with terminal value $d+1$, and write $j^{a\oplus}$ for its
$a$-fold iterate. In particular,
$j^{0\oplus}=j$ and $(d+1)^\oplus=d+1$.

\begin{example}\label{exam:betam}
Let
\[
  \overline w=(2,1,2,3,2,1)=(i_6,\ldots,i_1),\qquad
  v=s_2s_3s_1s_2.
\]
The rightmost reduced subexpression of $v^{-1}$ is
\[
  \overline v=(2,1,3,2)=(i_6,i_5,i_3,i_2),\qquad
  (p_1,p_2,p_3,p_4)=(2,3,5,6).
\]
For $k=p_4=6$, one has
\[
  f_{\min}(6)=1,\qquad f(6)=4,\qquad
  \gamma_4=2,\qquad b_4=1.
\]
\end{example}

\subsection{The initial coordinates for the completed word}

Choose and fix a completion
\[
  \dot v=(j_{n_0},\ldots,j_1),\qquad
  j_m=i_{p_m}\quad(1\leq m\leq d),\qquad
  s_{j_1}\cdots s_{j_{n_0}}=w_0.
\]
Let $\pi_v:\mathbb Z^{n_0}\to\mathbb Z^d$ be projection
onto the first $d$ coordinates. For an element $X$ with
a Lusztig parameter, write
\[
  \mathbf a^{\overline v}(X)
  :=\pi_v\bigl(\mathbf a^{\dot v}(X)\bigr).
\]
For each minor $D_k$, the selected and complementary positions
in Definition~\ref{def:extremal-weight-recursion} are those
associated with $u_k$ in the completed word $\dot v$.

\begin{proposition}\label{pro:amdotv}
For $1\leq k\leq r$ and $1\leq m\leq d$, one has
\begin{equation}\label{eq:amdotv3}
  a_m^{\dot v}(D_k)
  =a_{p_m}^{\dot w}(D_k)
  =
  \begin{cases}
    1,&p_m\leq k\text{ and }i_{p_m}=i_k,\\
    0,&\text{otherwise}.
  \end{cases}
\end{equation}
Consequently, $\mathbf a^{\overline v}(D_k)$ is independent
of the chosen completion of $\overline v$, and its support is
\[
  \{1\leq m\leq d\mid
      f_{\min}(k)\leq m\leq f(k),\ i_{p_m}=i_k\}.
\]
\end{proposition}

\begin{proof}
Put
\[
  s(k):=\#\{1\leq m\leq d\mid p_m\leq k\},
\]
and let $\mathcal Q_k$ be the set of positions of the
leftmost reduced subexpression of $u_k$ in $\dot v$.
By \cite[Proposition~5.27 and Lemma~5.29]{menard2022cluster},
with $w_{\mathrm M}=w^{-1}$ and $v_{\mathrm M}=v^{-1}$,
\[
  \mathcal Q_k\cap[1,d]=\{s(k)+1,\ldots,d\}.
\]
For $k=n_0$, this equality also holds, since
$\mathcal Q_k=\varnothing$ and $s(k)=d$.

Thus positions $1,\ldots,s(k)$ are complementary, and
\eqref{eq:xi} gives
\[
  \xi_m^{\dot v,k}
  =s_{j_1}\cdots s_{j_m}\varpi_{i_k}
  \qquad(0\leq m\leq s(k)).
\]
It follows that
\[
  \xi_{m-1}^{\dot v,k}-\xi_m^{\dot v,k}
  =\delta_{j_m,i_k}\beta_m^{\dot v}
  \qquad(1\leq m\leq s(k)).
\]
For $s(k)<m\leq d$, the position $m$ belongs to
$\mathcal Q_k$, so this difference is zero.
Since $j_m=i_{p_m}$ and $m\leq s(k)$ is equivalent to
$p_m\leq k$, Proposition~\ref{prop:extremal-minor-Lusztig-data}
gives the stated formula for $a_m^{\dot v}(D_k)$.
Its equality with $a_{p_m}^{\dot w}(D_k)$ follows from
\eqref{eq:interval-parameter}.
The remaining assertions follow from this formula and the
definitions of $f_{\min}$ and $f$.
\end{proof}

\subsection{M\'enard's mutation sequence}
\label{sec:mutationseq}

Let $\cS$ be the initial classical unipotent seed with
variables $D_k|_{q=1}$, and let $\cS^q$ be its quantum lift
with variables $D_k$, using the normalization fixed above.
Their additive counterparts are the summands $V_k$.
We apply the same sequence of mutations to these seeds,
keeping the original position labels throughout.
All mutations are computed in the full ambient seeds.

For the color $c=i_{p_m}$, list its occurrences in the
original word as
\[
  c_1<c_2<\cdots<c_{n_c}.
\]
Define the $m$-th block by
\begin{equation}\label{eq:mutation-block}
  \widetilde\mu_m
  :=\mu_{c_{n_c-\gamma_m}}\circ\cdots
       \circ\mu_{c_{b_m+1}},
\end{equation}
with the rightmost mutation performed first.
The block is the identity if $b_m+1>n_c-\gamma_m$.
Thus each block visits positions of one color in increasing
order. By \cite[Proposition~8.1]{menard2022cluster}, these
blocks agree with the recursive blocks used in the support
calculations of that paper.

Set
\[
  \widehat\mu_0=\mathrm{id},\qquad
  \widehat\mu_m=\widetilde\mu_m\circ\cdots\circ\widetilde\mu_1,
  \qquad
  \widetilde M=\widehat\mu_d,
\]
and write $\cS_m=\widehat\mu_m(\cS)$ and
$\cS_m^q=\widehat\mu_m(\cS^q)$.
In particular, $\cS_0=\cS$, and $\widetilde M=\mathrm{id}$
when $d=0$.

After block $m$, record as deleted the largest position
of color $i_{p_m}$ not already recorded as deleted.
The set recorded by this stage is
\[
  E_m=\{k\in[1,r]\mid
          k>(k_{\max})^{\alpha(k,m)-}\},
  \qquad E_0=\varnothing.
\]
These labels remain in the full ambient seed when later
mutations are computed. They are removed only from the
auxiliary cuts and, at the end, from the Richardson seed.
Put
\[
  E:=E_d,\qquad J:=[1,r]\setminus E.
\]

The auxiliary cuts also omit certain labels that will be
retained in the final seed. At an intermediate step $t$,
let $R_k^{(t)}$ be the additive summand with label $k$.
A label not recorded as deleted is called \emph{evicted}
if
\[
  R_k^{(t)}\in\cC^{\mathrm{add}}_{w,v}.
\]
By \cite[Proposition~5.35 and Definition~7.1]
{menard2022cluster}, this is equivalent to the vanishing
of the first $d$ coordinates of its additive
$\Delta$-vector for $\dot v$, and is precisely
M\'enard's eviction criterion.
The auxiliary cut is obtained by removing the recorded
and evicted labels, together with their incident arrows.
Eviction affects only this auxiliary cut; an evicted
label outside $E$ remains in the final Richardson seed.
The comparison in Section~\ref{sec:comparison} will
identify this additive criterion with the vanishing of
the truncated Lusztig parameter.

To obtain the final seed from $\widetilde M(\cS)$,
retain every original frozen vertex in $J$, freeze every
retained mutable vertex adjacent to a vertex of $E$,
and then delete $E$. Denote this seed by $\Sigma_{w,v}$,
its frozen index set by $J_{\mathrm{fr}}$, and its
variables by $(x_k)_{k\in J}$.
Apply the same operations to $\widetilde M(\cS^q)$
to obtain $\Sigma^q_{w,v}$.
This is a quantum seed: the deleted rows of the exchange
matrix vanish in every retained mutable column, so
restriction preserves compatibility with the quantum
commutation form.

We use the notation
\[
  \cA^{\mathrm{loc}}(\Sigma_{w,v})
  :=\cA(\Sigma_{w,v})
       [x_a^{-1}\mid a\in J_{\mathrm{fr}}],
\]
where $\cA(\Sigma_{w,v})$ is taken over $\mathbb C$
with noninvertible frozen variables.

\begin{example}\label{exam:mutation-A3}
For the data of Example~\ref{exam:betam}, the first selected
color is $2$, whose original positions are $2,4,6$.
Since $\gamma_1=1$ and $b_1=0$,
\[
  \widetilde\mu_1=\mu_4\circ\mu_2.
\]
The remaining three blocks are the identity. The final
deleted and retained sets are
\[
  E=\{3,4,5,6\},\qquad J=\{1,2\}.
\]
Thus $|J|=\ell(w)-\ell(v)=2$, and both retained vertices
are frozen.
\end{example}

\begin{theorem}[The classical Richardson seed]
\label{thm:classical-seed}
The seed $\Sigma_{w,v}$ is M\'enard's classical seed in
the conventions fixed above. Its frozen variables are
the cluster characters of all indecomposable
projective-injective objects of
$\cC^{\mathrm{add}}_{w,v}$, and
\begin{equation}\label{eq:classical-localized-seed}
  \cA^{\mathrm{loc}}(\Sigma_{w,v})
  \simeq\mathbb C[\Rwv].
\end{equation}
This isomorphism agrees with the geometric localization
of $\Iwv$ inside their common fraction field.
\end{theorem}

\begin{proof}
Apply \cite[Theorem~6.8 and Propositions~7.25--7.26]
{menard2022cluster} to $(v^{-1},w^{-1})$.
The surviving original frozen summands remain
projective-injective. Moreover, every retained vertex
isolated after deletion is either already frozen or
adjacent to a deleted vertex: otherwise it would be an
isolated ambient mutable vertex, whose additive exchange
sequence would have zero middle term.
Thus our freezing rule gives precisely M\'enard's seed.

Its frozen summands are the indecomposable
projective-injective objects of
$\cC^{\mathrm{add}}_{w,v}$, as in
\cite[Theorem~4.5]{leclerc2016cluster}.
Every cluster-tilting object contains all these objects,
so their cluster characters are the frozen variables in
every categorical seed.

The comparison with the weave cluster structure
\cite[Theorems~1.1 and~10.1]{casals2024cluster},
transported through \eqref{eq:left-right-richardson},
gives \eqref{eq:classical-localized-seed}.
Since these identifications preserve the specified
cluster-character functions, they induce the same
embedding into $\operatorname{Frac}(\Iwv)$ as the
geometric localization.
\end{proof}

M\'enard's theorem also shows that every final retained
additive summand belongs to $\cC^{\mathrm{add}}_{w,v}$.
Hence every such label is evicted from the auxiliary cut,
and the final auxiliary cut is empty. The corresponding
vanishing of quantum Lusztig parameters will be proved
in Section~\ref{sec:comparison}.

\begin{remark}
The original frozen vertices must be retained.
In type $A_1$, with $v=e$ and $w=s$, no vertex is deleted,
but the sole variable $x$ is frozen. Thus
\[
  \cA(\Sigma_{s,e})=\mathbb C[x],
  \qquad
  \mathbb C[\mathring{\mathcal B}_{e,s}]
     =\mathbb C[x,x^{-1}].
\]
Theorem~\ref{thm:classical-seed} concerns the latter
localization.
\end{remark}

\section{Additive categories and comparison of parameters}
\label{sec:comparison}

\subsection{Cluster characters and filtration parameters}

Retain the preprojective algebra $\Lambda$, the injective
modules $Q_i$, and the word conventions fixed above. Put
\[
  V_k=\soc_{(i_k,\ldots,i_1)}(Q_{i_k}),\qquad
  V_{\overline w}=\bigoplus_{k=1}^rV_k.
\]
Then $V_{\overline w}$ is a cluster-tilting object of
Leclerc's category $\mathcal C_{w^{-1}}$. The cluster
character satisfies
\[
  \varphi(V_k)=D_k|_{q=1},\qquad
  \varphi(M\oplus N)=\varphi(M)\varphi(N).
\]
We use this character to realize additive seeds as functions;
the ordinary Grothendieck group of the additive category
plays no role here.

For an $i$-box $[a,b]$, define
\[
  M^{\mathrm{add}}[a,b]:=V_b/V_{a^-},\qquad V_0:=0,
\]
using the natural inclusions between the $V_k$ of the same
color. In particular, put
$M_t^{\mathrm{add}}:=V_t/V_{t^-}$.
Every $R\in\mathcal C_{w^{-1}}$ admits an ordered GLS
filtration with layers that are direct sums of these modules.
Its multiplicities are uniquely determined; we denote their
vector by
\[
  \Delta_{\overline w}(R)
  =\bigl(\Delta_{\overline w,t}(R)\bigr)_{t=1}^r
  \in\mathbb Z_{\geq0}^r.
\]
More explicitly, it is characterized by
\[
  \bigl(\dim\operatorname{Hom}_\Lambda(V_s,R)\bigr)_{s=1}^r
  =
  \sum_{t=1}^r\Delta_{\overline w,t}(R)
  \bigl(\dim\operatorname{Hom}_\Lambda
       (V_s,M_t^{\mathrm{add}})\bigr)_{s=1}^r.
\]
The vectors on the right are linearly independent.
These facts follow from \cite[Section~10]{geiss2011kac};
see also \cite[Lemma~5.2, Propositions~5.4--5.5,
and Definition~5.7]{menard2022cluster}.
Only uniqueness of the multiplicity vector is needed.

The same construction applies to any reduced word.
For a completion to $w_0$, its domain is all finite-dimensional
$\Lambda$-modules. In particular, both
$\Delta_{\dot w}$ and $\Delta_{\dot v}$ are defined on every
module occurring below. For the initial summands,
\[
  \Delta_{\overline w,t}(V_k)
  =
  \begin{cases}
    1,&t\leq k\text{ and }i_t=i_k,\\
    0,&\text{otherwise}.
  \end{cases}
\]

\begin{proposition}[Initial truncated $\Delta$-vectors]
\label{pro_deltav}
For $1\leq k\leq r$ and $1\leq m\leq d$, one has
\[
  \Delta_{\dot v,m}(V_k)
  =\Delta_{\dot w,p_m}(V_k)
  =
  \begin{cases}
    1,&p_m\leq k\text{ and }i_{p_m}=i_k,\\
    0,&\text{otherwise}.
  \end{cases}
\]
\end{proposition}
\begin{proof}
This is \cite[Theorem~5.30]{menard2022cluster}, with
$v_{\mathrm M}=v^{-1}$ and $w_{\mathrm M}=w^{-1}$.
\end{proof}

For additive seeds, we use the Gabriel quiver of
$\operatorname{End}_\Lambda(T)^{\mathrm{op}}$.
Thus our arrows are the reversals of M\'enard's
\cite[Remark~2.3]{menard2022cluster}, and agree with the
exchange-quiver convention for $\cS$ and $\cS^q$.
The local statements used below are unchanged by this
simultaneous reversal.

Apply M\'enard's full mutation sequence to
$V_{\overline w}$ and to the quantum seed $\cS^q$.
At step $t$, denote the corresponding summands and quantum
variables by $R_k^{(t)}$ and $X_k^{(t)}$.
The initial identities and the exchange relations imply
\[
  \varphi(R_k^{(t)})=X_k^{(t)}|_{q=1}
  \qquad\text{for every }t\text{ and }k.
\]
As in Section~\ref{sec:mutationseq}, all position labels
remain in the full seeds, even after they have been
removed from an auxiliary cut.

\subsection{Comparison of truncated parameters}
\label{sec:deltavector}

Use the fixed completion $\dot v$ and projection $\pi_v$
from Section~\ref{sec:parameters}, and abbreviate
\[
  \mathbf d(R):=\pi_v\bigl(\Delta_{\dot v}(R)\bigr),
  \qquad
  \mathbf a_v(X):=\mathbf a^{\overline v}(X)
                 =\pi_v\bigl(\mathbf a^{\dot v}(X)\bigr).
\]
For $0\leq s\leq d$, put
$v_s=s_{j_1}\cdots s_{j_s}$ and let
$\mathcal F_s:=\mathcal C^{v_s^{-1}}$ be Leclerc's
torsion-free class. In particular, $v_0=e$ and $v_d=v$.
The prefix version of
\cite[Proposition~5.35]{menard2022cluster} gives
\begin{equation}\label{eq:prefix-torsionfree-criterion}
  R\in\mathcal F_s
  \quad\Longleftrightarrow\quad
  \Delta_{\dot v,a}(R)=0\quad(1\leq a\leq s)
\end{equation}
for every indecomposable rigid module under consideration.
Each $\mathcal F_s$ is closed under extensions and direct
summands. Since all ambient summands belong to
$\mathcal C_{w^{-1}}$, their eviction criterion from the
preceding section is precisely $\mathbf d(R)=0$.

If mutation at $k$ replaces $X_k$ by $X_k'$, then
Lemma~\ref{lem:leading-exchange} gives
\begin{equation}\label{eq:projected-exchange-candidates}
  \mathbf a_v(X_k)+\mathbf a_v(X_k')
  \in
  \left\{
    \sum_i[b_{ik}]_+\mathbf a_v(X_i),\
    \sum_i[-b_{ik}]_+\mathbf a_v(X_i)
  \right\}.
\end{equation}
Here the full PBW parameters are compared first, and the
selected value is then projected. In particular, this
statement does not take a coordinatewise maximum.

Let $Z_a$ be the summand at the label recorded as deleted
after block $a$. Such a summand remains in the full seed.
The next lemma explains why it does not enter later
exchange relations.

\begin{lemma}\label{lem:deleted-vertices-no-later-neighbours}
During block $m$, no label recorded as deleted in an
earlier block is adjacent in the full additive quiver
to a vertex being mutated.
\end{lemma}

\begin{proof}
Put $c=i_{p_a}$. Just before block $a$, the terminal
nondeleted label of this color is
$h=(h_{\max})^{\alpha(h,a-1)-}$.
The support formula in
\cite[Theorem~7.10(6)]{menard2022cluster} has lower endpoint
$f_{\min}(h)^{\alpha(h,a-1)\oplus}=a$ and upper endpoint
$f(h^{\alpha(h,a-1)+})=f(h_{\max})$.
Consequently,
\[
  \mathbf d(Z_a)
  =
  \sum_{\substack{a\leq b\leq d\\i_{p_b}=c}}\mathbf e_b,
\]
where $\mathbf e_b$ is the $b$-th standard coordinate vector.
The bounds in \eqref{eq:mutation-block} exclude $h$ from
block $a$ and every later block. Hence
\begin{equation}\label{eq:deleted-module-nonzero-coordinate}
  \Delta_{\dot v,a}(Z_a)=1
\end{equation}
continues to hold at every subsequent step.

Consider a mutation $R_k\leftrightarrow R_k'$ in block $m$.
By \cite[Theorem~7.10 and Lemma~7.14]{menard2022cluster},
both summands have zero coordinates in positions
$1,\ldots,m-1$; this also holds when $R_k'$ is evicted.
Thus both belong to $\mathcal F_{m-1}$.
Their two additive exchange sequences have the form
\[
  0\longrightarrow R_k\longrightarrow E_+
    \longrightarrow R_k'\longrightarrow0,\qquad
  0\longrightarrow R_k'\longrightarrow E_-
    \longrightarrow R_k\longrightarrow0.
\]
The middle terms are the sums of the neighbouring summands
in the two arrow directions, with multiplicities.
Closure under extensions and direct summands therefore
places every full-quiver neighbour of $k$ in
$\mathcal F_{m-1}$.
For $a<m$, however,
\eqref{eq:deleted-module-nonzero-coordinate} and
\eqref{eq:prefix-torsionfree-criterion} exclude $Z_a$
from this class. It cannot be a neighbour of $k$.
\end{proof}

For comparison with the additive cuts, define the
\emph{quantum auxiliary cut} by removing the same recorded
deleted labels and every other label with
$\mathbf a_v(X_k)=0$.

\begin{theorem}\label{thm:parameter-vanishing}
\label{prop:truncated-delta-lusztig-comparison}
Apply M\'enard's mutation sequence to $V_{\overline w}$
and to $\cS^q$. At every individual mutation step,
\begin{equation}\label{eq:truncated-delta-lusztig-comparison}
  \mathbf d(R_k^{(t)})=\mathbf a_v(X_k^{(t)})
  \qquad\text{for every label }k
\end{equation}
in the full seeds.
Consequently, the additive and quantum auxiliary cuts
have the same labels and exchange quivers in our conventions.
At the end of the sequence, every retained label $k\in J$
satisfies $\mathbf a_v(X_k)=0$.
\end{theorem}

\begin{proof}
The initial equality follows from
Propositions~\ref{pro_deltav} and~\ref{pro:amdotv}.
The full exchange quivers agree throughout the sequence,
since they agree initially and undergo the same mutations.

Assume the equality holds before mutation at $k$ in
block $m$, and suppress the step superscripts.
By Lemma~\ref{lem:deleted-vertices-no-later-neighbours},
previously deleted labels do not occur in either exchange
monomial at $k$. Evicted labels contribute zero to the
additive truncated sums and, by induction, to the quantum
truncated sums as well. We may therefore compute these
sums in the additive cut.

Call an arrow horizontal if its endpoints have the same
original color, and ordinary otherwise. By
\cite[Proposition~7.16]{menard2022cluster}, the horizontal
arrows at $k$ point in one direction and the ordinary
arrows in the other. Put
\[
  H_k:=\sum_{i\text{ joined horizontally to }k}\mathbf d(R_i),
  \qquad
  O_k:=\sum_{i\text{ joined ordinarily to }k}\mathbf d(R_i),
\]
counting arrow multiplicities. Equation
\eqref{eq:projected-exchange-candidates} and the induction
hypothesis give
\[
  \mathbf a_v(X_k')
  \in\{H_k-\mathbf d(R_k),\ O_k-\mathbf d(R_k)\}.
\]

The support formulas in
\cite[Theorem~7.10(6) and Lemma~7.14(6)]{menard2022cluster}
give $\mathbf d(R_k)_m=1$, whereas every ordinary
neighbour has $m$-th coordinate zero.
Thus
\[
  \bigl(O_k-\mathbf d(R_k)\bigr)_m=-1.
\]
Nonnegativity of Lusztig data excludes this alternative,
so
\[
  \mathbf a_v(X_k')=H_k-\mathbf d(R_k).
\]
On the additive side,
\cite[Lemma~7.17]{menard2022cluster} identifies the same
vector with $\mathbf d(R_k')$.
This proves the equality at the mutated label; all
other labels are unchanged.

Recording a deletion or evicting a label changes only
the auxiliary cut, so the induction applies to all full
seed labels. Blocks with no mutations require no extra
argument. Equality of the truncated vectors then
identifies the two auxiliary cuts.

Finally, Theorem~\ref{thm:classical-seed} places every
final retained additive summand in
$\cC^{\mathrm{add}}_{w,v}$.
Its vector $\mathbf d$ is therefore zero by
\eqref{eq:prefix-torsionfree-criterion}.
Equation~\eqref{eq:truncated-delta-lusztig-comparison}
gives the required vanishing for every $k\in J$.
\end{proof}

\section{Monoidal realization and localization}
\label{sec:categorification}

\subsection{Quiver Hecke categories and support}

Let $R=\bigoplus_{\nu\in Q_+}R(\nu)$ be a symmetric
quiver Hecke algebra of the same finite simply-laced type
over $\mathbb C$, with the standard geometric choice of
parameters. For example, for an orientation of the Dynkin
diagram one may take
\[
  Q_{ij}(u,v)=(v-u)^{m_{ij}}(u-v)^{m_{ji}}
  \qquad(i\neq j),
\]
where $m_{ij}$ is the number of arrows from $i$ to $j$.
We work with finite-dimensional graded modules.
Write $qM$ for grading shift, $M^*$ for graded duality,
and $M\circ N$ for convolution.

A simple module $M$ is \emph{real} if $M\circ M$ is simple.
Two simple modules \emph{strongly commute} if their
convolution is simple. We use self-dual representatives
of simple modules. For a pairwise strongly commuting
family of real self-dual modules, $\bigodot$ denotes the
self-dual grading normalization of their convolution.

Recall that $\mathbb A=\mathbb Z[q,q^{-1}]$, and write
\[
  \Kcl{\mathscr C}
  :=\mathbb C\otimes_{\mathbb A}K_0(\mathscr C),
  \qquad q\longmapsto1.
\]
For classical cluster algebras, monoidal categorification
requires an isomorphism with this specialized Grothendieck
ring under which cluster monomials are represented by
real simple objects. We will also realize every successive
seed mutation by exchange short exact sequences.

For a nonzero $R(\nu)$-module $M$, define
\begin{align*}
  \mathrm W(M)
  &:=\{\gamma\in Q_+\mid e(\gamma,\nu-\gamma)M\neq0\},\\
  \mathrm W^*(M)
  &:=\{\gamma\in Q_+\mid e(\nu-\gamma,\gamma)M\neq0\}.
\end{align*}
Here $e(\gamma,\delta)$ is the sum of the word idempotents
whose first block has root degree $\gamma$ and whose
second block has root degree $\delta$; it is zero unless
$\gamma,\delta\in Q_+$ and $\gamma+\delta=\nu$.
For modules with several root degrees, take the union of
the supports of their homogeneous summands.
Set $\mathrm W(0)=\mathrm W^*(0)=\varnothing$.
The positive root degree is $\nu$, while the usual
representation-theoretic weight is $-\nu$.

Define the full subcategories of $R\text{-gmod}$ by
\begin{align*}
  \mathscr C_w
  &:=\{M\mid \mathrm W(M)\subset
       \operatorname{span}_{\mathbb R_{\geq0}}
         (\Delta_+\cap w\Delta_-)\},\\
  \mathscr C_{*,v}
  &:=\{M\mid \mathrm W^*(M)\subset
       \operatorname{span}_{\mathbb R_{\geq0}}
         (\Delta_+\cap v\Delta_+)\},\\
  \mathscr C_{w,v}
  &:=\mathscr C_w\cap\mathscr C_{*,v}.
\end{align*}
These categories are closed under subquotients,
extensions, convolution, and grading shifts
\cite[Proposition~2.16]{kashiwara2018monoidal}.
The algebra $A_{w,v}$ retains the definition in
\eqref{eq:quantum-invariant-algebra}.

\begin{proposition}\label{prop:klr-identification}
Quiver Hecke categorification identifies
$K_0(R\text{-gmod})$ with
$A_q(\mathfrak n)_{\mathbb A}$, sending self-dual simple
modules to the upper global basis.
It restricts to a based identification
\[
  K_0(\mathscr C_{w,v})\simeq A_{w,v}
\]
and a classical isomorphism
\[
  \Phi_{\mathrm{cl}}:
  \Kcl{\mathscr C_{w,v}}\xrightarrow{\sim}\Iwv.
\]
The specialized classes of self-dual simple modules form
a $\mathbb C$-basis of $\Kcl{\mathscr C_{w,v}}$.
\end{proposition}
\begin{proof}
The ambient identification and its assertion about bases
are \cite[Theorems~2.1.2 and~2.1.4]{kang2018monoidal}.
Its restriction is
\cite[Theorem~2.20]{kashiwara2018monoidal}, and the
classical invariant algebra is identified in
\cite[Remarks~2.19 and~2.21(ii)]{kashiwara2018monoidal},
with the actions fixed in
\eqref{eq:quantum-invariant-algebra}.
Finally, the self-dual simple classes form an
$\mathbb A$-basis before specialization, and hence a
$\mathbb C$-basis after specialization.
\end{proof}

Write $\bfj=(j_1,\ldots,j_{n_0})$ for the forward reading
of the fixed completion $\dot v$.
The cuspidal multiplicities of a simple module $L$ for
this convex order equal the Lusztig parameters of its
upper global basis element
\cite{tingley2016mirkovic}.
We denote this common vector by
$\mathbf a^{\dot v}(L)=\mathbf a^{\bfj}(L)$ and use
$\mathbf a_v(L)$ for its truncation.

\begin{lemma}\label{lem:cuspidal-membership}
For a simple module $L\in\mathscr C_w$,
\[
  L\in\mathscr C_{w,v}
  \quad\Longleftrightarrow\quad
  \mathbf a_v(L)=0.
\]
For $d=0$, the coordinate condition is empty.
\end{lemma}
\begin{proof}
The roots attached to the first $d$ letters of $\bfj$
are precisely $\Delta_+\cap v\Delta_-$.
By \cite[Proposition~2.18]{kashiwara2018monoidal},
$L\in\mathscr C_{*,v}$ if and only if all its nonzero
cuspidal multiplicities occur after these positions.
The assertion follows.
\end{proof}

\begin{lemma}[Closure under convolution factors]
\label{lem:factor-closure}
Let $M,N$ be nonzero finite-dimensional graded modules.
If $M\circ N\in\mathscr C_{w,v}$, then
$M,N\in\mathscr C_{w,v}$.
\end{lemma}
\begin{proof}
By \cite[Lemma~2.2]{kashiwara2018monoidal},
\[
  \mathrm W(M\circ N)=\mathrm W(M)+\mathrm W(N),\qquad
  \mathrm W^*(M\circ N)=\mathrm W^*(M)+\mathrm W^*(N).
\]
Each factor support contains $0$, so it is contained in
the corresponding convolution support.
The defining cone conditions therefore hold for both factors.
\end{proof}

\subsection{The ambient seed and all retained mutations}

For every nonzero unipotent quantum minor
$D(u\lambda,v\lambda)$, let $M(u\lambda,v\lambda)$ be
the corresponding self-dual determinantial module.
In particular,
\[
  M^{\bfi}[a,b]
  :=M(w_{\leq b}\varpi_{i_a},w_{\leq a-1}\varpi_{i_a}),
  \qquad
  M^{\bfi}[a,a]=L_{\bfi}(\beta_a),
\]
where $L_{\bfi}(\beta_a)$ is the cuspidal module for the
convex order fixed by $\bfi$ and its completion.

The initial monoidal seed of $\mathscr C_w$ consists of
$M(w_{\leq k}\varpi_{i_k},\varpi_{i_k})$, $1\leq k\leq r$.
Every ambient seed reachable from it has a monoidal
realization in $\mathscr C_w$: its variables are real
simple modules, its normalized cluster monomials are
simple convolutions, and its mutations are realized by
exchange short exact sequences \cite{kang2018monoidal}.

Apply $\widetilde M$ to this ambient seed.
For each retained label $k\in J$, let $L_k$ be the
resulting self-dual simple module.
Theorem~\ref{thm:parameter-vanishing} and
Lemma~\ref{lem:cuspidal-membership} give
\begin{equation}\label{eq:seed-in-category}
  L_k\in\mathscr C_{w,v}\qquad(k\in J).
\end{equation}
Together with the restricted exchange and commutation
matrices, these modules define a monoidal seed
$\mathscr T_{w,v}$ lifting $\Sigma^q_{w,v}$.

We next verify that freezing and deletion remain
compatible with all subsequent retained mutations.

\begin{lemma}\label{lem:freezing-deletion}
Every mutation sequence in the mutable directions
$J_{\mathrm{ex}}:=J\setminus J_{\mathrm{fr}}$ of
$\Sigma^q_{w,v}$ is the restriction of the same mutation
sequence in the full ambient quantum seed.
At every step, the deleted rows of the ambient exchange
matrix vanish in all columns indexed by $J_{\mathrm{ex}}$.
\end{lemma}
\begin{proof}
The freezing rule gives $b_{ek}=0$ for $e\in E$ and
$k\in J_{\mathrm{ex}}$.
After mutation at $\ell\in J_{\mathrm{ex}}$,
$b'_{e\ell}=-b_{e\ell}=0$, and for $k\neq\ell$,
\[
  b'_{ek}
  =b_{ek}
   +[b_{e\ell}]_+[b_{\ell k}]_+
   -[-b_{e\ell}]_+[-b_{\ell k}]_+
  =0.
\]
This vanishing therefore persists.
The exchange monomials involve only retained variables,
and matrix mutation commutes with restriction.
The quantum toric frames consequently restrict as well.
Compatibility with the restricted commutation matrix
follows by omitting the zero deleted-row contributions
from the ambient compatibility identity.
\end{proof}

\begin{proposition}\label{prop:all-mutations}
Every mutation sequence in the mutable directions of
$\mathscr T_{w,v}$ is realized by monoidal seeds in
$\mathscr C_{w,v}$.
Consequently, every cluster monomial of
$\cA(\Sigma_{w,v})$, with nonnegative frozen exponents,
is the specialized class of a real simple module in
$\mathscr C_{w,v}$.
\end{proposition}
\begin{proof}
We induct on the number of mutations.
The initial modules belong to $\mathscr C_{w,v}$ by
\eqref{eq:seed-in-category}.
Suppose the assertion holds for the current seed and
mutate at $k\in J_{\mathrm{ex}}$.
Lemma~\ref{lem:freezing-deletion} identifies this with
the corresponding ambient mutation, whose exchange
monomials involve only retained variables.

Ambient monoidal categorification supplies a real simple
module $L_k'$ and an exchange sequence
\[
  0\longrightarrow q^{a}M_+
   \longrightarrow q^{b}(L_k\circ L_k')
   \longrightarrow q^{c}M_-\longrightarrow0,
\]
where $M_+$ and $M_-$ are normalized convolution
monomials in the current retained seed modules.
Both end terms belong to $\mathscr C_{w,v}$ by induction
and convolution closure. Extension closure places
$L_k\circ L_k'$ in the same category, and
Lemma~\ref{lem:factor-closure} gives
$L_k'\in\mathscr C_{w,v}$.
The opposite exchange sequence and all strong
commutation properties are inherited from the ambient seed.
This completes the induction.

Normalized convolution monomials in any resulting seed
are simple. Applying this to twice the exponent vector
shows that they are real.
Their specialized classes are the required cluster monomials.
\end{proof}

\begin{remark}\label{rem:no-admissibility-shortcut}
The construction realizes all successive mutations in
$\mathscr C_{w,v}$ before any equality between its
Grothendieck ring and a cluster algebra is used.
It also embeds the quantum cluster algebra with
noninvertible frozen variables into
$K_0(\mathscr C_{w,v})$ after the standard scalar
extension for half-integral grading normalizations.
No equality of these unlocalized quantum algebras
is asserted.
\end{remark}

\subsection{Determinantial and frozen localizations}
\label{subsec:monoidal-localization}

For $i\in I$, put
\[
  \mathsf D_i:=M(w\varpi_i,v\varpi_i).
\]
Its specialized class maps under $\Phi_{\mathrm{cl}}$
to the minor
\[
  d_i=\Delta_{v\varpi_i,w\varpi_i}
\]
from Proposition~\ref{prop:geometric-frozen-localization}.

We first account for the full-support hypothesis in
\cite[Section~4 and Theorem~4.4]{kashiwara2023localizations}.
If $w=e$, then $v=e$ and the assertions reduce to graded
vector spaces. Otherwise put $I'=\operatorname{supp}(w)$.
The subword criterion gives
$\operatorname{supp}(v)\subseteq I'$.
Every root degree in $\mathscr C_w$ belongs to
$Q_{I',+}$, so $\mathscr C_w$ and $\mathscr C_{w,v}$
identify with the corresponding categories for the
restricted Cartan datum on $I'$.
For $i\in I'$, restriction of the fundamental weight
gives the same determinantial module.
For $i\notin I'$, both $w$ and $v$ fix $\varpi_i$, and
$\mathsf D_i\simeq\mathbf1$.

The cited theorem thus makes $(\mathsf D_i)_{i\in I}$
a real commuting family of graded central objects of
$\mathscr C_{w,v}$. Here graded centrality means that
for homogeneous $M$ of positive root degree $\gamma$
there are natural isomorphisms
\[
  \mathsf D_i\circ M
  \simeq q^{\chi_i(\gamma)}M\circ\mathsf D_i,
\]
compatible with convolution, where $\chi_i$ is an
integer-valued linear function of the root degree.
Since determinantial modules are real and admit affinizations,
these isomorphisms imply that the $\mathsf D_i$ strongly
commute with every simple object of $\mathscr C_{w,v}$;
see \cite[Lemma~3.20 and Theorem~3.26]
{kashiwara2020localizations}.
We may therefore form
\[
  \widetilde{\mathscr C}_{w,v}
  :=\mathscr C_{w,v}
       [\mathsf D_i^{\circ-1}\mid i\in I],
  \qquad
  \mathsf Q:\mathscr C_{w,v}
       \longrightarrow\widetilde{\mathscr C}_{w,v}.
\]
By \cite[Proposition~2.13]{kashiwara2020localizations},
this is an abelian monoidal category with exact
convolution, and $\mathsf Q$ is exact.
Every localized object is, up to grading shift, of the
form $\mathsf Q(M)\circ\mathsf Q(D)^{\circ-1}$,
where $D$ is a convolution product of the $\mathsf D_i$.

For a simple $S\in\mathscr C_{w,v}$, each $D\circ S$
is simple. Clearing denominators in a nonzero morphism
to $\mathsf Q(S)$ represents it by a nonzero morphism
to such a simple module, which is an epimorphism.
Exactness of $\mathsf Q$ then shows that
$\mathsf Q(S)$ is simple.
It is nonzero because convolution with a nonzero
denominator cannot annihilate a nonzero module.
Applying $\mathsf Q$ to composition series and then
tensoring with invertible objects also shows that
$\widetilde{\mathscr C}_{w,v}$ has finite length.

\begin{lemma}\label{lem:frozen-denominator-comparison}
For $a\in J_{\mathrm{fr}}$, let $L_a$ be the frozen
module of $\mathscr T_{w,v}$ and put
$x_a:=\Phi_{\mathrm{cl}}([L_a])$.
There are integers $m_{ia}\geq0$ such that
\begin{equation}\label{eq:denominator-frozen-factorization}
  d_i=\prod_{a\in J_{\mathrm{fr}}}x_a^{m_{ia}},
  \qquad
  \mathsf D_i\simeq
    \bigodot_{a\in J_{\mathrm{fr}}}L_a^{\odot m_{ia}},
\end{equation}
and for every $a$ there is an $i$ with $m_{ia}>0$.
Consequently,
\[
  \Iwv[d_i^{-1}\mid i\in I]
  =
  \Iwv[x_a^{-1}\mid a\in J_{\mathrm{fr}}],
\]
and every $\mathsf Q(L_a)$ is invertible.
\end{lemma}
\begin{proof}
The common ambient mutation sequence identifies the
specialized monoidal frozen classes with the additive
cluster-character functions.
Proposition~\ref{prop:geometric-frozen-localization}
therefore gives the first factorization, the occurrence
of every frozen factor, and the equality of commutative
localizations.

The normalized product
$F_i:=\bigodot_aL_a^{\odot m_{ia}}$ is self-dual simple,
since the frozen modules are real and strongly commute.
Its specialized class has the same image $d_i$ as that
of the self-dual simple module $\mathsf D_i$.
By Proposition~\ref{prop:klr-identification}, these
specialized simple classes form a basis.
Hence $F_i\simeq\mathsf D_i$, proving the second
factorization.

Fix $a$ and choose $i$ with $m_{ia}>0$.
Moving one factor $L_a$ to either end gives
\[
  \mathsf D_i\simeq q^{u}L_a\circ Y,
  \qquad
  \mathsf D_i\simeq q^{v}Y\circ L_a
\]
for a frozen convolution monomial $Y$ and suitable
integers $u,v$.
Since $\mathsf Q(\mathsf D_i)$ is invertible, these
factorizations give a right inverse and a left inverse
of $\mathsf Q(L_a)$, respectively.
Thus $\mathsf Q(L_a)$ is invertible.
\end{proof}

\begin{theorem}\label{theo_categorification}
Assume finite simply-laced type and let $v\leq w$.
The localized category
$\widetilde{\mathscr C}_{w,v}$ is a monoidal
categorification of
$\cA^{\mathrm{loc}}(\Sigma_{w,v})\simeq\mathbb C[\Rwv]$.
More precisely, there is an isomorphism
\[
  \Kcl{\widetilde{\mathscr C}_{w,v}}
  \xrightarrow{\sim}\mathbb C[\Rwv]
\]
under which every cluster monomial with arbitrary
integer frozen exponents is represented by a real
simple object of $\widetilde{\mathscr C}_{w,v}$.
For nonnegative frozen exponents, such an object can
already be chosen in $\mathscr C_{w,v}$.
\end{theorem}

\begin{proof}
Proposition~\ref{prop:all-mutations} realizes every
retained mutation and every cluster monomial with
nonnegative frozen exponents in $\mathscr C_{w,v}$.

By Proposition~\ref{prop:klr-identification},
$K_0(\mathscr C_{w,v})$ embeds into the domain
$A_q(\mathfrak n)$.
Thus the classes $[\mathsf D_i]$ are regular.
Graded centrality makes them normal, with
\[
  [\mathsf D_i][M]
  =q^{\chi_i(\gamma)}[M][\mathsf D_i]
  \qquad
  \bigl(M\in\mathscr C_{w,v}\cap R(\gamma)\text{-gmod}\bigr).
\]
Their products therefore form an Ore denominator set,
allowing powers of $q$.
Fraction representatives and clearing common
denominators give
\[
  K_0(\widetilde{\mathscr C}_{w,v})
  \simeq
  K_0(\mathscr C_{w,v})
    \bigl[[\mathsf D_i]^{-1}\mid i\in I\bigr],
\]
by the same argument as in the proof of
\cite[Corollary~5.4]{kashiwara2020localizations}.
Indeed, fraction representatives give surjectivity;
preservation of simple objects and the isomorphism
criterion after clearing denominators give injectivity.

All commutation factors specialize to $1$.
The universal property of localization consequently gives
\[
  \Kcl{\widetilde{\mathscr C}_{w,v}}
  \simeq \Iwv[d_i^{-1}\mid i\in I]
  \simeq \mathbb C[\Rwv],
\]
where the second identification is
Proposition~\ref{prop:geometric-frozen-localization}.
By Lemma~\ref{lem:frozen-denominator-comparison}, this
is also localization at the frozen classes.
Theorem~\ref{thm:classical-seed} identifies the same
functions with
$\cA^{\mathrm{loc}}(\Sigma_{w,v})$.

The exact functor $\mathsf Q$ sends all exchange
sequences from Proposition~\ref{prop:all-mutations}
to exchange sequences in the localized category.
The seed objects remain nonzero real simple, and the
frozen objects are invertible.

Finally, write a cluster monomial with integer frozen
exponents as a monomial with nonnegative exponents
times an inverse frozen monomial.
The first factor is represented by a real simple object
$L\in\mathscr C_{w,v}$.
The second is represented by an invertible object that
commutes with $\mathsf Q(L)$ up to grading shift,
because all its factors belong to the same seed.
Tensoring with it preserves simplicity, and the same
commutation shows that the square of the resulting
object is simple.
This proves both the categorification statement and
the assertions about cluster monomials.
\end{proof}

\subsection{An explicit Richardson seed in type \texorpdfstring{$A_3$}{A3}}
\label{subsec:worked-richardson-example}

\begin{example}[A Richardson seed with mutable vertices]
\label{exam:mutable-richardson-A3}
Let $G=\mathrm{SL}_4$, $w=w_0$, and $v=s_2$.  Use the increasing
word
\[
 \bfi=(1,2,3,2,1,2),\qquad
 \overline w=(2,1,2,3,2,1).
\]
The selected position is $p_1=2$.  Thus $\gamma_1=1$,
$\beta_1=0$, and M\'enard's sequence is
$\widetilde M=\mu_4\circ\mu_2$.
Write a general element of $N$ as
\[
 n=\begin{pmatrix}
  1&a&b&c\\
  0&1&h&e\\
  0&0&1&f\\
  0&0&0&1
 \end{pmatrix},\qquad
 P=ahf-ae-bf+c,\quad F=bf-c.
\]
The initial flag minors, in position order, are
\[
 \left.(D_1,D_2,D_3,D_4,D_5,D_6)\right|_{q=1}
   =(a,ah-b,P,ae-c,c,be-ch),
\]
with frozen positions $\{3,5,6\}$.  The two ambient exchange
relations are the polynomial identities
\begin{align*}
 (ah-b)f&=(ae-c)+P,\\
 (ae-c)F&=af(be-ch)+Pc.
\end{align*}
Consequently, the full mutated classical cluster is
\[
 (a,f,P,F,c,be-ch).
\]
The deleted set is $E=\{6\}$, and the retained set is
$J=\{1,2,3,4,5\}$.  After these mutations, the only retained
mutable vertex adjacent to vertex $6$ is vertex $4$.
Thus the Richardson seed has
\[
 J_{\ex}=\{1,2\},\qquad J_{\fr}=\{3,4,5\},\qquad
 (x_1,x_2,x_3,x_4,x_5)=(a,f,P,F,c).
\]
With rows ordered by $1,\ldots,5$ and columns by $1,2$, its
extended exchange matrix, in our sign convention, is
\[
 B=\begin{pmatrix}
 0&0\\
 0&0\\
 -1&0\\
 1&1\\
 0&-1
 \end{pmatrix}.
\]

The principal part of $B$ is zero, so the cluster type is
$A_1\times A_1$. Let $X_k$ denote the corresponding self-dual
simple modules, with $x_k=\Phi_{\cl}([X_k])$.
For their Lusztig coordinates, choose the increasing completed word
$\bfj=(2,1,2,3,2,1)$, whose first letter represents $v$.
It is connected to $\bfi$ by
\[
 123212\longrightarrow123121
 \longrightarrow121321\longrightarrow212321.
\]
The usual rank-two Lusztig transition, together with the intervening
commutation, gives the initial first coordinates
$(0,1,0,1,0,1)$.  The two exchanges change this row successively to
\[
 (0,0,0,1,0,1)\quad\text{and}\quad(0,0,0,0,0,1).
\]
Indeed, the new coordinates at positions $2$ and $4$ are respectively
$\max\{1,0\}-1=0$ and $\max\{0+0+1,0+0\}-1=0$.
The classes $[X_k]$ are interval minors for $\bfj$ under the
based identification of Proposition~\ref{prop:klr-identification}.
Thus \eqref{eq:interval-parameter} checks their full parameters
directly:
\[
\begin{array}{c|c|c|c}
 k&x_k&[X_k]&\bfa^{\bfj}(X_k)\\\hline
 1&a&D^{\bfj}[3,3]&(0,0,1,0,0,0)\\
 2&f&D^{\bfj}[6,6]&(0,0,0,0,0,1)\\
 3&P&D^{\bfj}[4,4]&(0,0,0,1,0,0)\\
 4&F&D^{\bfj}[3,5]&(0,0,1,0,1,0)\\
 5&c&D^{\bfj}[2,6]&(0,1,0,0,0,1)
\end{array}
\]
Every retained first coordinate is zero, as required for
membership in $\mathscr C_{w_0,s_2}$.

There are genuine subsequent mutations.  For example, mutation at
vertex $1$ gives $x_1'=hf-e$. If $X_1'$ is its corresponding
self-dual simple module, then
\[
 a(hf-e)=P+F,\qquad
 \bfa^{\bfj}(X_1')=(0,0,0,0,1,0).
\]
Mutation at vertex $2$ similarly gives $x_2'=b$ through
$fb=F+c$.  The corresponding simple modules also have zero first coordinate.  Directly,
left multiplication by $x_2(t)$ replaces $(h,e)$ by
$(h+t,e+tf)$, so all the displayed retained variables and their
two exchange partners are $N(s_2)$-invariant.

Finally, generalized minors are ordinary subdeterminants here.
Their determinantial frozen factorizations are
\[
 \begin{aligned}
 d_1&=\Delta_{s_2\varpi_1,w_0\varpi_1}=c=x_5,\\
 d_2&=\Delta_{s_2\varpi_2,w_0\varpi_2}
      =\det\begin{pmatrix}b&c\\1&f\end{pmatrix}=F=x_4,\\
 d_3&=\Delta_{s_2\varpi_3,w_0\varpi_3}=P=x_3.
 \end{aligned}
\]
Thus every frozen variable is a determinantial denominator in this
example, and the corresponding self-dual modules satisfy
\[
 \mathsf D_1\simeq X_5,\qquad
 \mathsf D_2\simeq X_4,\qquad
 \mathsf D_3\simeq X_3.
\]
In fact, putting $g=hf-e$ gives
\[
 \mathcal I_{w_0,s_2}=\CC[a,b,c,f,g],\qquad
 \cA(\Sigma_{w_0,s_2})=\CC[a,b,c,f,g],
\]
since $F=bf-c$, $P=ag-F$, and the two exchanges produce $b$ and
$g$.  The zero principal part of $B$ makes these exchanges
independent, so the mutable variables are exactly $a,g,f,b$.
The invariant-ring equality also follows from the polynomial
change of coordinates $(h,e)\leftrightarrow(h,g)$, under which
$N(s_2)$ translates only $h$.  The localized coordinate ring is
therefore explicitly
\[
 \CC[\mathring{\mathcal B}_{s_2,w_0}]
   =\CC[a,b,c,f,g]
      [c^{-1},(bf-c)^{-1},(ag-bf+c)^{-1}].
\]
The unlocalized equality just established is specific to this
calculation and is not used in the general theorem.
\end{example}

\end{document}